\documentclass[a4paper,10pt]{article}

\usepackage{amsmath, amssymb, amsthm, mathtools}
\usepackage{cite}
\usepackage{color,soul}  
\usepackage{flafter}  
\usepackage{fullpage}
\usepackage{graphicx}
\usepackage[mathscr]{eucal} 
\usepackage[small]{caption}
\usepackage{subfigure}

\newif\ifupdatetikz
\updatetikzfalse 
\ifupdatetikz
	\usepackage{tikz}
	\usepackage{pgfplots}
	\pgfplotsset{
		compat=newest,
		tick label style={font=\scriptsize},
		label style={font=\footnotesize},
		legend style={font=\tiny}
	}
	\usepgfplotslibrary{external}
	\usetikzlibrary{shapes,arrows,positioning}
	\tikzset{cross/.style={cross out, draw=black, minimum size=2*(#1-\pgflinewidth), inner sep=0pt, outer sep=0pt},cross/.default={2pt}}
	\tikzstyle{dot} = [black,fill=black]
	\tikzstyle{red_dot} = [red,fill=red]
	\tikzstyle{axis} = [-latex',line width=1.25]
	\tikzstyle{stab_region} = [dot,fill=red!30,opacity=0.7]
	
\else
	\usepackage{tikzexternal}
\fi
\newcommand{\figname}[1]{\tikzsetnextfilename{#1}}

\newtheoremstyle{mystyle}
  {}
  {}
  {\itshape}
  {}
  {\bfseries}
  {.}
  { }
  {\thmname{#1}\thmnumber{ #2}\thmnote{ (#3)}}

\theoremstyle{mystyle}
\newtheorem{theorem}{Theorem}
\newtheorem{remark}[theorem]{Remark}

\let\originalleft\left
\let\originalright\right
\renewcommand{\left}{\mathopen{}\mathclose\bgroup\originalleft}
\renewcommand{\right}{\aftergroup\egroup\originalright}

\newcommand{\abs}[1]{\lvert#1\rvert}
\newcommand{\alphavec}{\boldsymbol{\alpha}}
\newcommand{\averageV}[1]{\left\langle #1 \right\rangle}
\newcommand{\Bcal}{\boldsymbol{\mathcal{B}}}
\newcommand{\betavec}{\boldsymbol{\beta}}
\newcommand{\C}{\mathbb{C}}			

\newcommand{\collOp}{\mathcal{Q}}
\newcommand{\collFreq}{\nu}
\newcommand{\collFreqvec}{\boldsymbol{\collFreq}}
\newcommand{\DF}{\mathbf{D}}
\newcommand{\Dt}{\Delta t}
\newcommand{\Dtime}{\mathrm{D}_t}
\newcommand{\Dxv}{\mathrm{D}_{\xGrid,\vGrid}}
\newcommand{\diag}[1]{\operatorname{diag}\left(#1\right)}
\newcommand{\diffCoef}{d}

\newcommand{\disk}{\mathcal{D}}

\newcommand{\dt}{\delta t}
\newcommand{\dx}{\Delta x}
\newcommand{\dy}{\Delta y}
\newcommand{\epsi}{\varepsilon}

\newcommand{\eye}{\mathbf{I}}

\newcommand{\fF}{\hat{\mathbf{F}}}
\newcommand{\fSDSys}{\mathbf{f}} 		

\newcommand{\hydroLimit}{\epsi \to 0}
\newcommand{\lambdavec}{\boldsymbol{\lambda}}

\newcommand{\Lw}{L_{\omega}}
\newcommand{\lw}{l_{\omega}}

\newcommand{\Ma}{\mathcal{M}}						

\newcommand{\MaxBGK}[2]{\Ma_{#1}\left(#2\right)} 	
\newcommand{\MaxBGKOneD}{\MaxBGK{v}{f}} 				
\newcommand{\MaxBGKMultiD}{\MaxBGK{\v}{f}} 			
\newcommand{\MaxGlobParams}[3]{\Ma^{#1,#2,#3}_{v}} 		
\newcommand{\MaxBGKLin}{\MaxBGK{\text{\scriptsize{lin}},v}{f}} 		
\newcommand{\MaxBGKLinDiscVSys}{\MaxBGK{\text{\scriptsize{lin}},j}{f}}	
\newcommand{\MaxBGKLinSD}{\MaxBGK{\text{\scriptsize{lin}},\vGrid}{\fSDSys}}	
\newcommand{\MaxBGKLinMultiD}{\MaxBGK{\text{\scriptsize{lin}},\v}{f}} 		

\newcommand{\MF}{\mathbf{M}}

\newcommand{\omegabar}{\bar{\omega}}
\newcommand{\PF}{\mathbf{P}}
\newcommand{\R}{\mathbb{R}}				

\newcommand{\SzeroF}{\mathbf{S}_{0}}
\newcommand{\set}[2]{\left(#1\right)_{#2=1}^{\uppercase{#2}}}

\newcommand{\sigmazerovec}{\boldsymbol{\sigma}_0}

\newcommand{\Tinf}{T^{\infty}}

\newcommand{\V}{V}

\newcommand{\Vmat}{\mathbf{V}}
\renewcommand{\v}{\mathbf{v}}
\newcommand{\vMacroOneD}{\bar{v}}
\newcommand{\vMacroMultiD}{\mathbf{\bar{v}}}
\newcommand{\vGrid}{\boldsymbol{v}}
\newcommand{\vbarinf}{\bar{v}^{\infty}}

\newcommand{\x}{\mathbf{x}}
\newcommand{\xGrid}{\boldsymbol{x}}

\newcommand{\zetavec}{\boldsymbol{\zeta}}

\newcommand{\Cnm}[2]{\C^{#1 \times #2}}
\newcommand{\Rnm}[2]{\R^{#1 \times #2}}

\begin{document}

\title{Telescopic projective integration for kinetic equations with multiple relaxation times}

\author{
	Ward Melis \thanks{Department of Computer Science, K.U. Leuven, Celestijnenlaan 200A, 3001 Leuven, Belgium ({\tt ward.melis@cs.kuleuven.be}).} 
\and 
	Giovanni Samaey \thanks{Department of Computer Science, K.U. Leuven, Celestijnenlaan 200A, 3001 Leuven, Belgium ({\tt giovanni.samaey@cs.kuleuven.be}).}}

\maketitle

\begin{abstract}	
	We study a general, high-order, fully explicit numerical method for simulating kinetic equations with a BGK-type collision model with multiple relaxation times. In that case, the problem is stiff and its spectrum consists of multiple separated eigenvalue clusters. Projective integration methods are explicit integration schemes that first take a few small (inner) steps with a simple, explicit method, after which the solution is extrapolated forward in time over a large (outer) time step. These are very efficient schemes, provided there are only two clusters of eigenvalues. Telescopic projective integration methods generalize the idea of projective integration methods by constructing a hierarchy of projective levels. Here, we show how telescopic projective integration methods can be used to efficiently integrate kinetic equations with multiple relaxation times. We show that the required number of projective levels depends on the number of clusters, which in turn depends on the stiffness of the BGK source term. The size of the outer level time step only depends on the slowest time scale present in the model and is independent of the stiffness of the problem. We discuss stability and illustrate the approach with simulations in one and two spatial dimensions.
\end{abstract}

\section{Introduction} \label{sec:intro}
The Boltzmann equation forms the cornerstone of the kinetic theory of rarefied gases. In a scalar $D$-dimensional setting without external forcing, this equation portrays the evolution of the one-particle distribution function $f(\x,\v,t)$ in phase space $(\x,\v)$ at time $t$, and takes the following general form \cite{cercignani1988boltzmann}:
\begin{equation} \label{eq:kin_eq_general}
	\partial_t f + \v\cdot\nabla_{\x} f = \collOp(f),
\end{equation}
in which $\x = \set{x^d}{d} \in \R^D$ and $\v = \set{v^d}{d} \in \R^D$ denote the positions and velocities of the particles, respectively.
The left hand side of equation \eqref{eq:kin_eq_general} represents the transport of particles with velocity $\v$, whereas the collision operator $\collOp$ on the right hand side describes velocity changes as a consequence of collisions between particles. The collision operator proposed by Boltzmann is the most general and considers collisions between any two particles \cite{cercignani1988boltzmann}. However, due to the resulting quadratic cost, it is very expensive to discretize. To reduce the computational cost, the full Boltzmann collision operator is typically replaced by the well-known BGK model \cite{Bhatnagar1954}, in which collisions are modeled as a linear relaxation towards a local Maxwellian equilibrium distribution function $\MaxBGKMultiD$:
\begin{equation} \label{eq:kin_eq_bgk} 
	\partial_t f + \v\cdot\nabla_{\x} f = \frac{\collFreq}{\epsi}(\MaxBGKMultiD - f),
\end{equation}
where $\collFreq(\x,t) \in \R^{+}$ is the collision frequency and $0 < \epsi \ll 1$ is a positive small-scale parameter that determines the relaxation time scale. The specific expression of $\collFreq(\x,t)$ depends on the dimension of velocity space and the type of microscopic collisions considered \cite{cercignani1988boltzmann}. For instance, for Maxwellian molecules, in a one-dimensional velocity space one simply uses $\collFreq=1$, whereas in two dimensions, one takes $\collFreq(\x,t) = \rho(\x,t)$ with $\rho(\x,t)$ the density of particles, see, e.g., \cite{Struchtrup2005}.

The difficulty of numerically integrating equations of the form \eqref{eq:kin_eq_bgk} follows from the stiffness present on the right hand side, for which appropriate numerical methods need to be selected. 
There is currently a large research effort in the design of algorithms that are uniformly stable in $\epsi$ and approach a scheme for the limiting equation when $\epsi$ tends to 0; such schemes are called asymptotic-preserving in the sense of Jin \cite{Jin1999}. We briefly review here some achievements, and refer to the cited references for more details. In \cite{Jin1998,Jin2000a,Klar1998}, separating the distribution function $f$ into its odd and even parts in the velocity variable results in a coupled system of transport equations where the stiffness appears only in the source term, allowing to use a time-splitting technique \cite{Strang1968} with implicit treatment of the source term; see also related work in \cite{Jin1999,Klar1999,Klar1999a}. Implicit-explicit (IMEX) schemes are an extensively studied technique to tackle this kind of problems \cite{Ascher1997,Filbet2011} (and references therein). Recent results in this setting were obtained by Dimarco et al. to deal with nonlinear collision kernels \cite{Dimarco2013asymptotic}, and an extension to hyperbolic systems in a diffusive limit is given in \cite{Boscarino2013}. A different point of view based on well-balanced methods was introduced by Gosse and Toscani \cite{Gosse2003,Gosse2004}, see also \cite{Buet2006,Buet2007}. Discontinuous Galerkin schemes have also been developed \cite{Adams2001,Guermond2010,Larsen1989,Lowrie2000,McClarren2008}, as well as regularization methods \cite{Haack2008,Hauck2009}. When the collision operator allows for an explicit computation, an explicit scheme can be obtained subject to a classical diffusion CFL condition by splitting $f$ into its mean value and the first-order fluctuations in a Chapman-Enskog expansion form \cite{Godillon-Lafitte2005}. Also, closure by moments, e.g. \cite{Coulombel2005}, can lead to reduced systems for which time-splitting provides new classes of schemes \cite{Carrillo2008}, see \cite{Lemou2008,Minerbo1978,Pomraning1991,Struchtrup2005} for more complete references on moment methods in general. Alternatively, a micro-macro decomposition based on a Chapman-Enskog expansion has been proposed \cite{Lemou2008}, leading to a system of transport equations that allows to design a semi-implicit scheme without time splitting. An non-local procedure based on the quadrature of kernels obtained through pseudo-differential calculus was proposed in \cite{Besse2010}. Finally, we refer to \cite{Dimarco2014} for a clear survey on numerical methods for kinetic equations.

A robust and fully explicit method, which allows for time integration of (two-scale) stiff systems with arbitrary order of accuracy in time, is projective integration. 
Projective integration was proposed in \cite{Gear2003projective} for stiff systems of ordinary differential equations with a clear gap in their eigenvalue spectrum. In such stiff problems, the fast modes, corresponding to the Jacobian eigenvalues with large negative real parts, decay quickly, whereas the slow modes correspond to eigenvalues of smaller magnitude and are the solution components of practical interest. Projective integration allows a stable yet explicit integration of such problems by first taking a few small (inner) steps using a step size $\dt$ with a simple, explicit method, until the transients corresponding to the fast modes have died out, and subsequently projecting (extrapolating) the solution forward in time over a large (outer) time step of size ${\Dt > \dt}$. In \cite{Lafitte2012}, projective integration was analyzed for kinetic equations with a diffusive scaling. An arbitrary order version, based on Runge-Kutta methods, has been proposed recently in \cite{Lejon2016}, where it was also analyzed for kinetic equations with an advection-diffusion limit \cite{Melis2015}. Alternative approaches to obtain a higher-order projective integration scheme have been proposed in \cite{Lee2007,Rico-Martinez}. 
These methods fit within recent research efforts on numerical methods for multiscale simulation \cite{E2003a,E2007,Kevrekidis2003,Kevrekidis2009}.

Projective integration methods work best whenever the problem's spectrum consists of two eigenvalue clusters (one corresponding to the fast and the other to the slow eigenvalues) with a clear spectral gap between them. When the outer time step is much larger than the inner time step, $\Dt \gg \dt$, the stability domain of projective integration methods consists of two circle-like stability regions \cite{Gear2003projective}. In that case, the projective integration method parameters can be tuned such that (i) all fast eigenvalues of the naive time discretization of the kinetic equation fall into its first stability region, and (ii) its dominant stability region contains all dominant eigenvalues of this naive time discretization. However, when more relaxation time scales are present, the spectrum will contain more than two eigenvalue clusters. Since the projective integration method possesses only two stability regions at most, the only way that stability of the method can be guaranteed is by choosing the method parameters such that its stability region does not split up into two parts and contains both the fastest and slowest eigenvalues. In that case, the projective integration method is called $[0,1]$-stable \cite{Gear2003telescopic}. Unfortunately, this requirement completely destroys much of the potential speed-up of the method, since this imposes a severe condition on the maximum possible value of the projective time step, thus defeating its purpose \cite{Gear2003projective}.

To integrate problems with multiple eigenvalue clusters, telescopic projective integration (TPI) methods, which are presented in \cite{Gear2003telescopic}, can be used. In these methods, the outer integrator step of the classical projective integration method is seen as the inner integrator of yet another outer integrator on a coarser level. By repeating this idea, TPI methods construct a hierarchy of projective levels in which each outer integrator step on a certain level serves as an inner integrator step one level higher. TPI methods can remedy the aforementioned difficulty of multiple spectral gaps in two distinct ways resulting in different criteria for selecting the method parameters: (i) they can be designed such that the method is always $[0,1]$-stable with a greater speed-up than classical projective integration, or alternatively (ii) they can be set up such that there is one stability region around every eigenvalue cluster. The latter is discussed in more detail in \cite{Gear2001}. To conclude this literature review, we also refer to \cite{Eriksson2004,Sommeijer1990,Vandekerckhove2007} for related approaches.

In this paper, we will construct telescopic projective integration methods of arbitrary order of accuracy in time to integrate kinetic equations of the form \eqref{eq:kin_eq_bgk}. As an intermediate step, which is interesting in its own right, we consider the relaxation time of the collisions to vary in space only. This variation is embodied in a relaxation profile function denoted by $\omega(\x) \in \R$. In that case, we end up with the following multiple relaxation times kinetic equation:
\begin{equation} \label{eq:kin_eq_omega}
	\partial_t f + \v\cdot\nabla_{\x} f = \frac{\omega(\x)}{\epsi}(\MaxBGKMultiD - f).
\end{equation}
The introduction of a relaxation profile function $\omega(\x)$ leads to a time-invariant spectrum that, in general, comprises multiple eigenvalue clusters separated by spectral gaps. The relaxation profile in equation \eqref{eq:kin_eq_omega} can, for instance, be understood as the mathematical description of a composite material in which each material has its own properties, which naturally leads to differing collisional relaxation times.

The remainder of this paper is structured as follows. In section \ref{sec:model}, we introduce our mathematical setup in more detail, discuss the choice of the Maxwellian function $\MaxBGKMultiD$ in 1D and 2D and compute the spectrum of both equations \eqref{eq:kin_eq_bgk} and \eqref{eq:kin_eq_omega} in 1D. In section \ref{sec:telescopic_proj_int}, we describe the telescopic projective integration method that will be used to integrate these kinetic equations. In section \ref{sec:num_props}, we determine the TPI method parameters for solving equation \eqref{eq:kin_eq_omega} based on its spectrum and extend the TPI construction procedure to kinetic equations of the form \eqref{eq:kin_eq_bgk}. We report numerical results in section \ref{sec:results}, where we illustrate the method for equation \eqref{eq:kin_eq_bgk} in one and two space dimensions. We conclude in section \ref{sec:conclusions}.

\section{Model problem} \label{sec:model}

\subsection{Kinetic equation and linearization} \label{subsec:kinetic}
In this work, we are interested in the BGK-type kinetic equation \eqref{eq:kin_eq_bgk} describing the evolution of a nonnegative one-particle distribution function $f(\x,\v,t) \in \R$, in which the particle positions and velocities are denoted by $\x \in \R^D$ and $\v \in \V \subset \R^D$, respectively. We are specifically interested in the cases $D=1$ and $D=2$. The right hand side of equation \eqref{eq:kin_eq_bgk} represents the BGK collision operator \cite{Bhatnagar1954}, modeling linear relaxation of the distribution function $f(\x,\v,t)$ to a local Maxwellian distribution $\MaxBGKMultiD \in \R$.
We also introduce the position density ${\rho(\x,t)=\averageV{f(\x,\v,t)}}$, obtained via averaging over the measured symmetric velocity space $(\V,\mu)$, 
\begin{equation}\label{eq:bracket} 
	\rho := \averageV{f} = \int_{\V} f d\mu(\v).
\end{equation} 
The local Maxwellian equilibrium function $\MaxBGKMultiD$ in the BGK model corresponds to the equilibrium distribution function of the full Boltzmann operator and is given by:
\begin{equation} \label{eq:maxwellian} 
	\MaxBGKMultiD = \frac{\rho}{(2\pi T)^{D/2}} \exp\left({-\frac{\abs{\v-\vMacroMultiD}^2}{2T}}\right),
\end{equation}
in which the density $\rho(\x,t) \in \R^{+}$, the mean velocity $\vMacroMultiD(\x,t) = \set{\vMacroOneD^d(\x,t)}{d} \in \R^D$ and the temperature $T(\x,t) \in \R^{+}$ are obtained as the moments of the distribution function $f$,
\begin{equation} \label{eq:moments_of_f} 
	\rho = \averageV{f}, \qquad 
	\vMacroOneD^d = \frac{1}{\rho}\averageV{v^d f}, \qquad
	T = \frac{1}{D\rho}\averageV{\abs{\v-\vMacroMultiD}^2 f}.
\end{equation}

\paragraph{One-dimensional case.}
To simplify the analysis in this work, we focus on the one-dimensional case $(D=1)$. In that case, equation \eqref{eq:kin_eq_bgk} reads:
\begin{equation} \label{eq:kin_eq_bgk_1d}
	\partial_t f + v\partial_x f = \frac{\collFreq}{\epsi}(\MaxBGKOneD - f),
\end{equation}
in which the particle positions and velocities are denoted by $x \in \R$ and $v \in \V \subset \R$, respectively.

According to the BGK model, particles interact with a collision frequency $\collFreq(x,t)$ that depends on the dimension of the velocity space. 
Since we are interested in studying the numerical difficulties that arise when dealing with multiple relaxation times, we choose $\collFreq(x,t)=\rho(x,t)$.  While this choice only improves the modeling accuracy of the BGK approximation in 2D, compared to simply choosing $\nu=1$, the numerical difficulties that are associated with the appearance of the additional relaxation times are the same (but easier to analyze) in 1D as in 2D.  

To facilitate the calculations of the spectrum (section \ref{subsec:spectrum}) and the numerical simulations (section \ref{sec:results}), we introduce a linearized Maxwellian distribution, denoted by $\MaxBGKLin$, obtained by linearizing the local Maxwellian distribution in equation \eqref{eq:maxwellian} around the following distribution,
\begin{equation} \label{eq:distribution}
	\MaxGlobParams{\rho}{\vbarinf}{\Tinf} = \frac{\rho}{\sqrt{2\pi\Tinf}}\exp\left({-\frac{|v - \vbarinf|^2}{2\Tinf}}\right),
\end{equation}
where $\rho$ is given in equation \eqref{eq:moments_of_f}, and $\vbarinf$ and $\Tinf$ are parameters of the distribution. If we evaluate this linearization for $\vbarinf = 0$, $\Tinf = 1$, and for a constant mean velocity $\vMacroOneD = 1$ and background temperature $T=1$ this leads to the following linearized Maxwellian:
\begin{equation} \label{eq:maxwellian_artificial} 
	\MaxBGKLin = \frac{\rho(1 + v)}{\sqrt{2\pi}}\exp\left({-\frac{v^2}{2}}\right),
\end{equation}
from which we define the velocity measure as: 
\begin{equation} \label{eq:V_measure}
	\mu(v) = \frac{1}{\sqrt{2\pi}}\exp\left(-\frac{v^2}{2}\right).
\end{equation}
We remark that this linearized Maxwellian falls into the class of Maxwellian distributions that was studied in \cite{Bouchut1999} in the setting of kinetic equations as relaxation models for hyperbolic conservation laws, see also \cite{Aregba-Driollet2000}. For a projective integration method in this context, we refer to \cite{Melis2015}.
In \cite{Melis2015}, it is shown that in the hydrodynamic limit, $\hydroLimit$, and on long time scales, equation \eqref{eq:kin_eq_bgk_1d} with the linearized Maxwellian \eqref{eq:maxwellian_artificial} tends to the dynamics of the linear advection equation:
\begin{equation} \label{eq:lin_adv}
	\partial_t \rho + \partial_x \rho = \epsi\partial_x(\diffCoef\partial_x \rho),
\end{equation}
where the right hand side of equation \eqref{eq:lin_adv} contains a small diffusive term with diffusion coefficient $\epsi\diffCoef$. Since the exact solution of the linear advection equation is known, the linearized Maxwellian case provides a means of assessing the accuracy of the proposed numerical technique.


\paragraph{Two-dimensional case.} In two space dimensions, equation \eqref{eq:kin_eq_bgk} is written as:
\begin{equation} \label{eq:kin_eq_bgk_2d}
	\partial_t f + v^x\partial_x f + v^y\partial_y f = \frac{\collFreq}{\epsi}(\MaxBGKMultiD - f),
\end{equation}
in which $\x = (x,y) \in \R^2$ and ${\v = (v^x,v^y) \in \V \subset \R^2}$ denote the particle positions and velocities, respectively. In this case, based on the one-dimensional linearized Maxwellian \eqref{eq:maxwellian_artificial}, we postulate the following Maxwellian distribution in 2D:
\begin{equation} \label{eq:maxwellian_artificial_2d} 
	\MaxBGKLinMultiD = \frac{\rho(1 + v^x)(1 + v^y)}{2\pi} \exp\left({-\frac{\abs{\v}^2}{2}}\right),
\end{equation}
from which we derive the velocity measure as:
\begin{equation} \label{eq:V_measure_2d}
	\mu(\v) = \frac{1}{2\pi}\exp\left(-\frac{\abs{\v}^2}{2}\right).
\end{equation}
For the Maxwellian in equation \eqref{eq:maxwellian_artificial_2d}, the dynamics of equation \eqref{eq:kin_eq_bgk_2d} in the hydrodynamic limit, $\hydroLimit$, and on long time scales now tends to the dynamics of the two-dimensional linear advection equation:
\begin{equation} \label{eq:lin_adv_2d}
	\partial_t \rho + \partial_x \rho + \partial_y \rho = O(\epsi),
\end{equation}
see also \cite{Melis2015}. Therefore, we can also compare with the exact solution in 2D.

In what follows, we will always assume that the velocity space is discrete, symmetric and of the form
\begin{equation} 
	\V := \set{\v_j}{j}, \qquad d\mu(\v)=\sum_{j=1}^J w_j \delta(\v-\v_j)d\v,
\end{equation}
where the chosen velocities satisfy $\v_{j} = -\v_{J-j+1}$ and $w_j$ represent the corresponding weights for which we have $\sum_{j=1}^J w_j = 1$. These discrete velocities $\set{\v_j}{j}$ and weights $\set{w_j}{j}$ are derived from the measures given in \eqref{eq:V_measure} and \eqref{eq:V_measure_2d} as the nodes and weights of the corresponding Gauss-Hermite quadrature. In this case, equations~\eqref{eq:maxwellian_artificial} or \eqref{eq:maxwellian_artificial_2d} break up into a system of $J$ coupled partial differential equations,
\begin{equation} \label{eq:system_part_diff_eq} 
	\partial_t f_j + \v_j\cdot\nabla_\x f_j = \frac{\collFreq}{\epsi}(\MaxBGKLinDiscVSys - f_j), \qquad  1 \le j \le J,
\end{equation}
in which $f_j(\x,t) \equiv f(\x,\v_j,t)$.

\subsection{Spectrum of the linearized kinetic equation} \label{subsec:spectrum}
To analyze stability of the telescopic projective integration technique, we need to investigate in more detail the spectrum of (a spatial discretization of) the one-dimensional BGK-type kinetic equation \eqref{eq:kin_eq_bgk_1d}, together with the linearized Maxwellian $\MaxBGKLin$ given in equation \eqref{eq:maxwellian_artificial}. To that end, we first discretize the system of equations \eqref{eq:system_part_diff_eq} on a uniform, constant in time, periodic spatial mesh with spacing $\dx$, consisting of $I$ mesh points $x_i=i\dx$, $1 \le i \le I$, with $I\dx=1$. After discretizing in space, we obtain the following semi-discrete system of ordinary differential equations:
\begin{equation}\label{eq:semidiscrete}
	\dot{\fSDSys} = \Dtime(\fSDSys),  \qquad \Dtime(\fSDSys):=- \Dxv(\fSDSys) + \frac{\collFreqvec}{\epsi}(\MaxBGKLinSD - \fSDSys),
\end{equation}
in which $\fSDSys$ and $\collFreqvec$ are vectors of length $I \times J$ resulting from the discretization in space and velocity, and $\Dxv(\cdot)$ represents a suitable discretization of the convective derivative $v\partial_x$ (e.g., upwind differences) where $\xGrid = \set{x_i}{i}$ and $\vGrid = \set{v_j}{j}$ denote the discrete grids in space and velocity, respectively.

As an intermediate step, in section \ref{subsubsec:spectrum_omega} we calculate the spectrum of system \eqref{eq:semidiscrete} when considering a time-invariant collision frequency $\collFreq(x,t) = \omega(x)$. Afterwards, we extend the obtained results to the case $\collFreq(x,t) = \rho(x,t)$ in section \ref{subsubsec:spectrum_rho}.

\subsubsection{Time-invariant collision frequency} \label{subsubsec:spectrum_omega}
We begin by deriving the spectrum of the semi-discrete system \eqref{eq:semidiscrete} using $\collFreq(x,t)=\omega(x)$ with a constant relaxation profile function $\omega(x) = \omegabar \in \R^{+}$. System \eqref{eq:semidiscrete} then becomes:
\begin{equation}\label{eq:semidiscrete_omega}
	\dot{\fSDSys} = \Dtime(\fSDSys),  \qquad \Dtime(\fSDSys):=- \Dxv(\fSDSys) + \frac{\omegabar}{\epsi}(\MaxBGKLinSD - \fSDSys).
\end{equation}
We assume that $\omegabar$ is bounded below by: $0 < \epsi \ll \omega_{\textrm{min}} \le \omegabar$, with $\omega_{\textrm{min}}$ independent of $\epsi$ such that there is a clear spectral gap. We transform the semi-discrete system of equations \eqref{eq:semidiscrete_omega} to the (spatial) Fourier domain yielding:
\begin{equation}\label{eq:fourier_semidiscrete}
	\partial_t \fF(\zeta_i) = \Bcal\;\fF(\zeta_i), \qquad
	\Bcal = \dfrac{\omegabar}{\epsi}\left(\MF\PF - \eye + \frac{\epsi}{\omegabar}\DF \right),
\end{equation}
in which $\fF \in \C^J$, $\Bcal$, $\DF \in \Cnm{J}{J}$, $\MF$, $\PF \in \Rnm{J}{J}$, and $\eye$ represents the identity matrix of dimension $J$. 
In equation \eqref{eq:fourier_semidiscrete}, the matrix $\DF$ represents the (diagonal) Fourier matrix of the spatial discretization chosen for the convection part, which depends on the Fourier mode $\zeta_i=2\pi i\dx$, $\PF$ is the Fourier matrix of the averaging of $f$ over the discrete velocity space $\V$, and the matrix $\MF$ corresponds to the Fourier transform of the linearized Maxwellian in equation \eqref{eq:maxwellian_artificial},
\[ \MF = \eye + \Vmat, \]
with $\Vmat$ the diagonal matrix given by $\diag{\vGrid}$, and using the velocity measure defined in \eqref{eq:V_measure}.

Since the velocity space is symmetric, we have the following property on the diagonal elements of the matrix $\DF$: 
\begin{equation}\label{eq:ass2} 
	D_{J-j+1}=\bar{D_j}, \qquad 1\le j \le J/2.
\end{equation}
Moreover, we write, from now on, 
\[ D_j = \alpha_j + \imath\beta_j, \]
in which $\alphavec = \set{\alpha_j}{j}$ and $\betavec = \set{\beta_j}{j}$ depend on the spatial discretization technique, the Fourier mode $\zeta_i$ and the chosen velocity grid $\vGrid$.
The following theorem is a corollary to \cite[Theorem 4.1]{Melis2015}.
\begin{theorem} \label{thm_spec}
	Under the above assumptions, the spectrum of matrix ${\Bcal = \dfrac{\omegabar}{\epsi}\left(\MF\PF - \eye + \dfrac{\epsi}{\omegabar}\DF \right)}$ satisfies
	\begin{equation} \label{eq:spectrum_B} 
	\mathrm{Sp}(\Bcal)
		\subset
	\left\{\disk\left(-\dfrac{\omegabar}{\epsi}, \max_{1 \le j \le J}\left(\sqrt{\alpha_j^2 + \beta_j^2}\right)\right)
		\cup
	\left\{\lambda^{(1)}(\omegabar)\right\}\right\},
	\end{equation}
	in which $\disk(c,r)$ denotes the disk with center $(c,0)$ and radius $r$. 
	The dominant eigenvalue $\lambda^{(1)}(\omegabar)$ is simple and can be expanded as 
	
	\begin{align}
		\Re\left\{\lambda^{(1)}(\omegabar)\right\} &= \averageV{\alphavec} + \left(\averageV{\alphavec^2} - \averageV{\alphavec}^2 - \averageV{\betavec^2} + \averageV{\betavec\vGrid}^2\right)\frac{\epsi}{\omegabar} + O\left(\frac{\epsi^2}{\omegabar^2}\right), \label{eq:dom_eig_real} \\
		\Im\left\{\lambda^{(1)}(\omegabar)\right\} &=  \averageV{\betavec\vGrid} + O\left(\frac{\epsi^2}{\omegabar^2}\right). \label{eq:dom_eig_imag}
	\end{align}
	
\end{theorem}

\begin{figure}[t]
	\begin{center}
		\figname{spectrum_equation_omega}
%
%
\newcommand{\figWidth}{14cm} 
\newcommand{\figHeight}{3cm} 
\newcommand{\subfigWidth}{\figWidth/4}
\newcommand{\subfigSpacingRight}{\subfigWidth/2} 
\begin{tikzpicture}
\begin{axis}[%
width=\subfigWidth,
height=\figHeight,
at={({2*(\subfigWidth+\subfigSpacingRight)},0)},
scale only axis,
xmin=-16000,
xmax=0,
xtick={-15000,-10000,-5000,0},
every outer y axis line/.append style={white},
every y tick label/.append style={font=\color{white}},
ymin=-1000,
ymax=1000,
ytick={\empty},
axis background/.style={fill=white},
axis x line*=bottom,
axis y line*=left
]
\addplot [color=blue,only marks,mark=x,mark options={solid},forget plot]
  table[]{tikz/data/spectrum_eq_omega-1.tsv};
\addplot [color=red,solid,forget plot,each nth point=10]
  table[]{tikz/data/fast_eig_zones-3.tsv};
\addplot [color=red,solid,forget plot,each nth point=10]
  table[]{tikz/data/fast_eig_zones-4.tsv};
\end{axis}

\draw (2*\subfigWidth+\subfigSpacingRight,0) -- ({2*(\subfigWidth+\subfigSpacingRight)},0);

\begin{axis}[%
width=\subfigWidth,
height=\figHeight,
at={(\subfigWidth+\subfigSpacingRight,0)},
scale only axis,
xmin=-207000,
xmax=-193000,
xtick={-205000, -200000, -195000},
xlabel={$\text{Re(}\lambda\text{)}$},
every outer y axis line/.append style={white},
every y tick label/.append style={font=\color{white}},
ymin=-1000,
ymax=1000,
ytick={\empty},
axis background/.style={fill=white},
axis x line*=bottom,
axis y line*=left
]
\addplot [color=blue,only marks,mark=x,mark options={solid},forget plot]
  table[]{tikz/data/spectrum_eq_omega-3.tsv};
\addplot [color=red,solid,forget plot,each nth point=10]
  table[]{tikz/data/fast_eig_zones-2.tsv};
\end{axis}

\draw (\subfigWidth,0) -- (\subfigWidth+\subfigSpacingRight,0);

\begin{axis}[%
width=\subfigWidth,
height=\figHeight,
at={(0,0)},
scale only axis,
xmin=-1007000,
xmax=-993000,
xtick={-1005000, -1000000,  -995000},
x tick label style={/pgf/number format/.cd, fixed, fixed zerofill, precision=3, /tikz/.cd},
ymin=-1000,
ymax=1000,
/pgf/number format/1000 sep={},
ytick = {-1000, -500, ..., 1000},
ylabel={$\text{Im(}\lambda\text{)}$},
axis background/.style={fill=white},
axis x line*=bottom,
axis y line*=left
]
\addplot [color=blue,only marks,mark=x,mark options={solid},forget plot]
  table[]{tikz/data/spectrum_eq_omega-6.tsv};
\addplot [color=red,solid,forget plot,each nth point=10]
  table[]{tikz/data/fast_eig_zones-1.tsv};
\end{axis}

\draw (0,\figHeight/2) -- (\figWidth,\figHeight/2);

\draw (\subfigWidth+\subfigSpacingRight/2+0.1cm,-0.2cm) -- ++(110:0.4cm);
\draw (\subfigWidth+\subfigSpacingRight/2+0.25cm,-0.2cm) -- ++(110:0.4cm);

\draw (\subfigWidth+\subfigSpacingRight/2+0.1cm,\figHeight/2-0.2cm) -- ++(110:0.4cm);
\draw (\subfigWidth+\subfigSpacingRight/2+0.25cm,\figHeight/2-0.2cm) -- ++(110:0.4cm);

\draw (2*\subfigWidth+3*\subfigSpacingRight/2+0cm,-0.2cm) -- ++(110:0.4cm);
\draw (2*\subfigWidth+3*\subfigSpacingRight/2+0.15cm,-0.2cm) -- ++(110:0.4cm);

\draw (2*\subfigWidth+3*\subfigSpacingRight/2+0cm,\figHeight/2-0.2cm) -- ++(110:0.4cm);
\draw (2*\subfigWidth+3*\subfigSpacingRight/2+0.15cm,\figHeight/2-0.2cm) -- ++(110:0.4cm);

\end{tikzpicture}%
	\end{center}
	\vspace{-0.4cm}\caption{\label{fig:spectrum_eq_omega} Spectrum of the semi-discretized system \eqref{eq:semidiscrete} with $\collFreq(x,t)=\omega(x)$ and $\omega(x)$ piecewise constant containing 4 $\omega$-values verifying the formal result in equation \eqref{eq:spectrum_Btilde}. For clarity, the scaling of the real axis is adapted to clearly visualize the different parts of the spectrum. }
\end{figure}
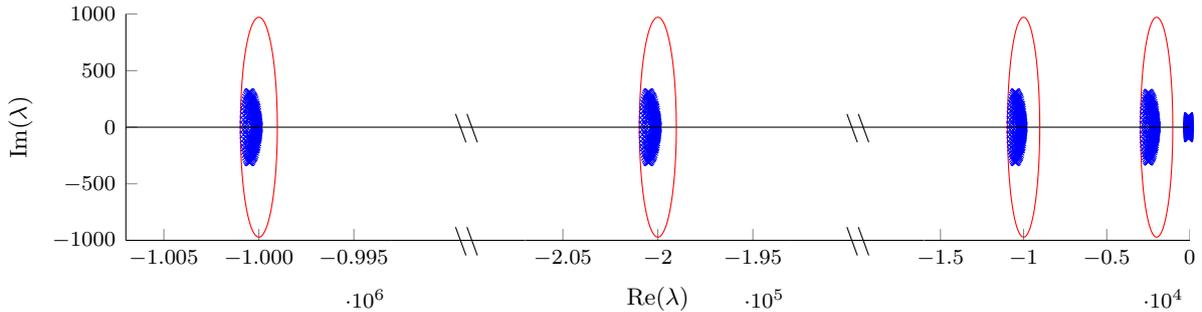

When the relaxation profile function $\omega(x)$ is piecewise constant over the spatial domain consisting of $\Lw$ constant values $\omegabar_{\lw},\;\lw=1,...,\Lw$, the above expressions in the spatial Fourier domain rapidly become very difficult for in this case the equations contain convolutions. However, we can still qualitatively identify the spectrum of the (formal) amplification matrix, which we denote by $\tilde\Bcal$, by performing numerical experiments. These experiments suggest that the spectrum of $\tilde\Bcal$ consists of (i) a combination of the $\Lw$ fast spectra obtained when considering constant $\omega(x)=\omegabar_{\lw}$ for all $x$ and for every $\lw=1,...,\Lw$, and (ii) the dominant eigenvalues $\lambda^{(1)}(\omegabar_1)$ given by \eqref{eq:dom_eig_real}-\eqref{eq:dom_eig_imag}, in which $\omegabar_1$ is considered to be the largest $\omega$-level. Therefore, using the spectrum of $\Bcal$ in \eqref{eq:spectrum_B}, we formally write the following conjecture on the spectrum of $\tilde\Bcal$:
\begin{equation} \label{eq:spectrum_Btilde}
\textrm{Sp}(\tilde\Bcal)
	\subset
\left\{\bigcup_{\lw=1}^{\Lw}\disk\left(-\dfrac{\omegabar_{\lw}}{\epsi}, \max_{1 \le j \le J}\left(\sqrt{\alpha_j^2 + \beta_j^2}\right)\right)
	\cup
\left\{\lambda^{(1)}(\omegabar_1)\right\}\right\}.
\end{equation}
Let us numerically illustrate this result. To that end, we set $\epsi=10^{-6}$ and discretize velocity space using $J=10$ velocity components corresponding to the Gauss-Hermite quadrature nodes for integration with respect to the measure in equation \eqref{eq:V_measure}. We consider $x\in[0, 1]$, apply periodic boundary conditions and use the upwind scheme of order $1$ with grid spacing $\dx=0.01$ as spatial discretization technique. The relaxation profile $\omega(x)$ is chosen as a piecewise constant function containing 4 well separated $\omega$-levels $\{1, 0.2, 0.01, 0.002\}$ that are distributed in zones of equal length over the spatial domain. The resulting spectrum is shown in figure \ref{fig:spectrum_eq_omega}. From this, we observe 4 fast eigenvalue clusters corresponding to the 4 $\omega$-values which are positioned around $\omegabar_{\lw}/\epsi$, $\lw=1,...,4$ and 1 slow cluster in the neighborhood of 0. The red disks $\disk(-\omegabar_{\lw}/\epsi,R_f)$ bound the fast eigenvalues where the radius $R_f$ is calculated as the maximal radius in equation \eqref{eq:spectrum_Btilde}:
\[ R_f = \max_{\zetavec,\vGrid}\left(\sqrt{\alphavec^2 + \betavec^2}\right), \]
with $\zetavec = \set{\zeta_i}{i}$. In this illustration, we have $R_f = 971.89$.

\subsubsection{Time-varying collision frequency} \label{subsubsec:spectrum_rho}
Moving on to the time-varying case $\collFreq(x,t) = \rho(x,t)$, the same reasoning as in the previous section can be used to obtain the eigenvalue spectrum of system \eqref{eq:semidiscrete}. 

When $\rho(x,t)$ is a continuous function, the spectrum would be continuously spread along the negative real axis in the interval $[-\max_{x} \rho(x,t)/\epsi,0]$. Notice that the size of this continuous spectral interval is time-dependent due to the time dependency of its left end point.
However, by discretizing in space we obtain a discrete spectrum where each value $\rho(x_i,t)$ can be seen as an $\omega$-level of the previous section. Extending the same reasoning used to derive equation \eqref{eq:spectrum_Btilde} we now find the following result:
\begin{equation} \label{eq:spectrum_Btilde_rho}
\textrm{Sp}(\tilde\Bcal)
	\subset
\left\{\bigcup_{i=1}^{I}\disk\left(-\dfrac{\rho(x_i,t)}{\epsi}, \max_{1 \le j \le J}\left(\sqrt{\alpha_j^2 + \beta_j^2}\right)\right)
	\cup
\left\{\lambda^{(1)}(\bar\rho(t))\right\}\right\},
\end{equation}
in which $\tilde\Bcal$ now refers to the amplification matrix of system \eqref{eq:semidiscrete} in the Fourier domain using ${\collFreq(x,t) = \rho(x,t)}$. We remark that the centers of the fast eigenvalue disks now depend on the space-discretized particle density $\rho(x_i,t)$. In that respect, a special case arises when $\rho(x_i,t) \to 0$ for some space-time points $(x_i,t)$ as this implies a transition from a fast cluster to an additional slow cluster.

\begin{figure}[t]
	\begin{center}
		\figname{spectrum_equation_rho}
%
%
\newcommand{\figHeight}{3cm} 
\newcommand{\figSpacingTop}{1.5cm} 
\begin{tikzpicture}


\begin{axis}[%
width=0.55\textwidth,
height=\figHeight,
at={(0,-(\figHeight+\figSpacingTop))},
scale only axis,
xmin=-1000000,
xmax=1000,
xlabel={$\text{Re(}\lambda\text{)}$},
ymin=-400,
ymax=400,
ytick = {-400, -200, ..., 400},
ylabel={$\text{Im(}\lambda\text{)}$},
axis background/.style={fill=white},
axis x line*=bottom,
axis y line*=left
]
\draw[solid, draw=red,thick] (axis cs:-10000,-400) rectangle (axis cs:0,400);
\addplot [color=black,solid,forget plot]
  table[]{tikz/data/spectrum_eq_rho-3.tsv};
\addplot [color=blue,only marks,mark=x,mark options={solid},forget plot]
  table[]{tikz/data/spectrum_eq_rho-4.tsv};
\end{axis}

\begin{axis}[%
width=0.25\textwidth,
height=\figHeight,
at={(0.65\textwidth,-(\figHeight+\figSpacingTop))},
scale only axis,
xmin=-10000,
xmax=0,
xlabel={$\text{Re(}\lambda\text{)}$},
ymin=-400,
ymax=400,
ytick = {-400, -200, ..., 400},
ylabel={$\text{Im(}\lambda\text{)}$},
axis background/.style={fill=white}
]
\addplot [color=black,solid,forget plot]
  table[]{tikz/data/spectrum_eq_rho-1.tsv};
\addplot [color=blue,only marks,mark=x,mark options={solid},forget plot]
  table[]{tikz/data/spectrum_eq_rho-2.tsv};
\end{axis}

\end{tikzpicture}%
	\end{center}
	\vspace{-0.4cm}\caption{\label{fig:spectrum_eq_rho} Spectrum of the semi-discretized system \eqref{eq:semidiscrete} with $\collFreq(x,t)=\rho(x,t)$ depicted for the initial density $\rho(\xGrid,0)$ in equation \eqref{eq:spectrum_rho_rho_init}. The left plot shows the global spectrum whereas the right plot provides a closer look of the slow part of the spectrum within the red rectangle. This confirms the formal result in equation \eqref{eq:spectrum_Btilde_rho}. }
\end{figure}
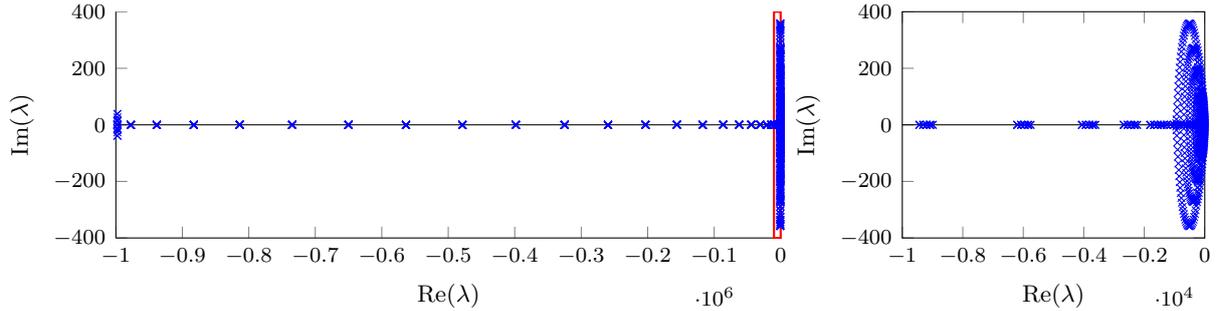

A numerical experiment is performed in which we use the same parameters as in the previous section. Since, in this case, the spectrum evolves with time we only plot the result for an initial density given by:
\begin{equation} \label{eq:spectrum_rho_rho_init}
	\rho(\xGrid,0) = \exp(-100(\xGrid-0.5)^2).
\end{equation}
The spectrum is shown in figure \ref{fig:spectrum_eq_rho}. We observe that the spectrum is indeed spread along the negative real axis in the interval $[-\max_{\xGrid} \rho(\xGrid,0)/\epsi,0]$ with $\max_{\xGrid} \rho(\xGrid,0) = 1$ for equation \eqref{eq:spectrum_rho_rho_init}. In addition, we find a number of extra slow clusters in the right hand side plot of figure \ref{fig:spectrum_eq_rho} corresponding to values of $\rho(\xGrid,0)$ that are sufficiently close to 0.

\section{Telescopic projective integration} \label{sec:telescopic_proj_int}
In this section, we construct a fully explicit, asymptotic-preserving, arbitrary order time integration method for the stiff semi-discretized system \eqref{eq:semidiscrete} containing in general more than two distinct time scales. The asymptotic-preserving property \cite{Jin1999} implies that, in the limit when $\epsi$ tends to zero, an $\epsi$-independent time step constraint, of the form $\Dt = O(\dx)$, can be used, similar to the hyperbolic CFL constraint for the limiting equation \eqref{eq:lin_adv}. To achieve this and overcome the difficulties mentioned in the introduction, we will use an extension of the projective integration method that can handle multiple time scales, entitled telescopic projective integration (TPI) \cite{Gear2003telescopic}.

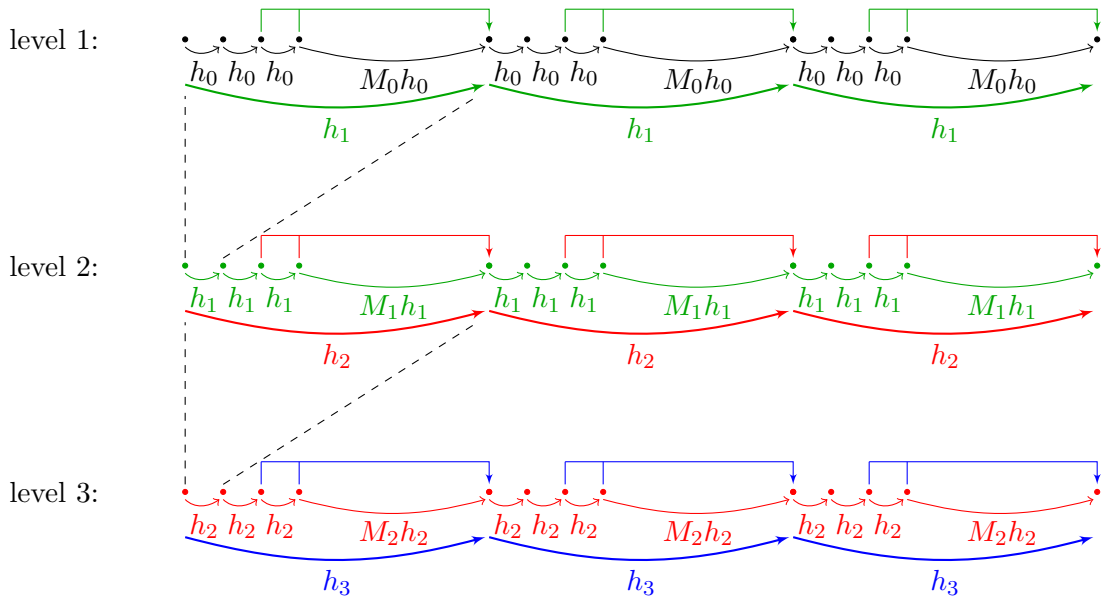
\begin{figure}[t]
	\begin{center}
		\figname{TPI_sketch}
		\newcommand{\TPIdots}[5] 
{
	\draw[dot,#4] (#1,#5) circle (.2ex); \draw[dot,#4] (#1+#3,#5) circle (.2ex); \draw[dot,#4] (#1+2*#3,#5) circle (.2ex); \draw[dot,#4] (#1+3*#3,#5) circle (.2ex);
}
\newcommand{\TPIdotsarrow}[6] 
{
	\TPIdots{#1}{#2}{#3}{#4}{#6};
	\draw [#5] (#1+2*#3,#6+0.1) -- (#1+2*#3,#6+0.4); \draw [#5] (#1+3*#3,#6+0.1) -- (#1+3*#3,#6+0.4); \draw [-latex',#5] (#1+2*#3,#6+0.4) -- (#2,#6+0.4) -- (#2,#6+0.1);
	\draw[dot,#4] (#2,#6) circle (.2ex);
}
\newcommand{\TPIdotsarrows}[7]
{
	\TPIdotsarrow{#1}{#2}{#3}{#4}{#5}{#7};
	\draw [->,#4] (#1,#7-0.1) to [out=-45,in=-145] (#1+#3-0.05,#7-0.1); \node [below,#4] at (#1+#3/2,#7-0.15) {$h_#6$};
	\draw [->,#4] (#1+#3,#7-0.1) to [out=-45,in=-145] (#1+2*#3-0.05,#7-0.1); \node [below,#4] at (#1+3*#3/2,#7-0.15) {$h_#6$};
	\draw [->,#4] (#1+2*#3,#7-0.1) to [out=-45,in=-145] (#1+3*#3-0.05,#7-0.1); \node [below,#4] at (#1+5*#3/2,#7-0.15) {$h_#6$};
	\draw [->,#4] (#1+3*#3,#7-0.1) to [out=-15,in=-165] (#2-0.05,#7-0.1); \node [below,#4] at (#1/2+3*#3/2+#2/2,#7-0.25) {$M_#6h_#6$};
	
	\edef\indexcurrenttimestep{#6};  \pgfmathparse{\indexcurrenttimestep+1};  \edef\indexnexttimestep{\pgfmathresult};
	\pgfmathparse{#2-#1};  \edef\length{\pgfmathresult};
	\draw[-latex',thick,#5] (#1,#7-0.6) to [out=-15,in=-165] (#2-0.05,#7-0.6); \node[below,#5] at (#1+\length/2,#7-0.9) {$h_{\pgfmathprintnumber\indexnexttimestep}$};
}

\begin{tikzpicture}[scale=1]
	\newcommand{\LevelHeight}{3}
	\newcommand{\OuterDistance}{4} 
	\newcommand{\InnerDistance}{0.5} 
	\node[left] at (-1,2*\LevelHeight) {level 1:};
	\TPIdotsarrows{0}{\OuterDistance}{\InnerDistance}{black}{green!65!black}{0}{2*\LevelHeight};
	\TPIdotsarrows{\OuterDistance}{2*\OuterDistance}{\InnerDistance}{black}{green!65!black}{0}{2*\LevelHeight};
	\TPIdotsarrows{2*\OuterDistance}{3*\OuterDistance}{\InnerDistance}{black}{green!65!black}{0}{2*\LevelHeight};

	\node[left] at (-1,\LevelHeight) {level 2:};
	\draw [black,dashed] (0,\LevelHeight+0.1) -- (0,2*\LevelHeight-0.75); \draw [black,dashed] (\InnerDistance,\LevelHeight+0.1) -- (\OuterDistance-0.15,2*\LevelHeight-0.75); 
	\TPIdotsarrows{0}{\OuterDistance}{\InnerDistance}{green!65!black}{red}{1}{\LevelHeight};
	\TPIdotsarrows{\OuterDistance}{2*\OuterDistance}{\InnerDistance}{green!65!black}{red}{1}{\LevelHeight};
	\TPIdotsarrows{2*\OuterDistance}{3*\OuterDistance}{\InnerDistance}{green!65!black}{red}{1}{\LevelHeight};

	\node[left] at (-1,0) {level 3:};
	\draw [black,dashed] (0,0+0.1) -- (0,\LevelHeight-0.75); \draw [black,dashed] (\InnerDistance,0.1) -- (\OuterDistance-0.15,\LevelHeight-0.75); 
	\TPIdotsarrows{0}{\OuterDistance}{\InnerDistance}{red}{blue}{2}{0};
	\TPIdotsarrows{\OuterDistance}{2*\OuterDistance}{\InnerDistance}{red}{blue}{2}{0};
	\TPIdotsarrows{2*\OuterDistance}{3*\OuterDistance}{\InnerDistance}{red}{blue}{2}{0};
\end{tikzpicture}
	\end{center}
	\vspace{-0.4cm}\caption{\label{fig:tpi_sketch} A level-3 telescopic projective integration method drawn for three outermost time steps $h_3$. The dots correspond to the different time points at which the numerical solution is calculated. We used $K=2$ constant on all levels. The time step and projective step size of each level $\ell=0,...,2$ are denoted by $h_\ell$ and $M_\ell$, respectively. }
\end{figure}

Projective integration is a time integration method that allows a stable yet explicit integration of two-scale stiff problems by first taking a few small (inner) steps with a simple, explicit method, until the transients corresponding to the fast modes have died out, and subsequently projecting (extrapolating) the solution forward in time over a large (outer) time step \cite{Gear2003projective}. 
Telescopic projective integration builds on this idea by employing a number of such projective integrator levels, which, starting from a base (\emph{innermost}) integrator, are wrapped around the previous level integrator \cite{Gear2003telescopic}. In this way, a hierarchy of projective integrators is formed in which each level fulfills both an inner and outer integrator role (except, of course, for the innermost and outermost level which only serve as an inner and outer integrator, respectively). This generalizes the idea of projective integration, which contains only one projective level wrapped around an inner integrator. On that account, in the TPI framework, the projective integration method is called a level-1 TPI method. The idea of a level-3 TPI method is sketched in figure \ref{fig:tpi_sketch}.

The different level integrators can in principle be selected independently from each other, but in general one selects a first order explicit scheme 
(e.g., the forward Euler scheme) for all but the outermost integrator level, whose order is chosen to meet certain accuracy requirements dictated by the problem.

The remainder of this section is structured as follows. In sections \ref{subsec:innermost} and \ref{subsec:outers} we present the innermost integrator and the projective outer levels, respectively. We summarize the TPI method parameters in section \ref{subsec:tpi_params} and discuss stability of the method in section \ref{subsec:tpi_stability}.

\subsection{Innermost integrator} \label{subsec:innermost}
We intend to integrate the semi-discrete system of equations \eqref{eq:semidiscrete} using a uniform time mesh with time step $h_0$, i.e., $t^k=kh_0$. The innermost integrator of the TPI method is chosen to be an explicit scheme, for which we use the following shorthand notation:
\begin{equation} \label{eq:timestepper_innermost}
	\fSDSys^{k+1} = S_{0}(\fSDSys^{k}),\qquad k = 0, 1, \ldots,
\end{equation}
in which $S_0$ denotes the time stepper with corresponding time step $h_0$. The forward Euler (FE) method immediately comes to mind, for which equation \eqref{eq:timestepper_innermost} is written as:
\begin{equation} \label{eq:FE_scheme} 
	\fSDSys^{k+1} = \fSDSys^k + h_0\Dtime(\fSDSys^k). 
\end{equation}
The purpose of the innermost integrator is to capture the fastest components in the numerical solution of system \eqref{eq:semidiscrete} and to sufficiently damp these out. We only require the innermost integrator to be stable for these components. Nevertheless, higher-order extensions of equation \eqref{eq:FE_scheme} such as the Runge-Kutta methods of order 2 and 4 are possible. However, as observed in \cite{Melis2015}, these higher-order methods bring forth severe stability restrictions on the projective integrator wrapped around the innermost integrator, in particular on the number of innermost time steps that is required per first level projective time step. Furthermore, the discretization error of the full TPI method will be dominated by the error of the outermost integrator. Consequently, we will not consider higher-order methods as innermost integrators in this work.

\subsection{Projective (outer) levels} \label{subsec:outers}
The telescopic projective integration method employs in general $L$ nested projective levels that are constructed around the innermost integrator as its fundamental building block. We now provide the scheme of the method in a framework similar to that of classical projective integration. Alternatively, the scheme can also be formulated recursively, see \cite{Gear2003telescopic}.

To keep track of the time instant at which the numerical solution is computed throughout the TPI method and at the same time desiring a compact notation, in what follows we employ superscript triplets of the form $(\ell,n,k_\ell)$ where $\ell$ denotes the integrator level ranging from $0$ (innermost) to $L-1$, $n$ represents the index of the current outermost integrator time $t^n=nh_L$, and $k_\ell$ corresponds to the iteration index of the integrator on level $\ell$. The numerical time on each level $\ell$ is then defined as:
\begin{equation} \label{eq:tpi_time}
	t^{\ell,n,k_\ell} = nh_L + \sum_{\ell'=\ell}^{L-1} k_{\ell'}h_{\ell'}.
\end{equation}
Notice that this time requires the iteration indices $k_{\ell'}$ of all outer integrators of a certain level $\ell$. Therefore, it incorporates a memory that keeps up with the current time instants at which the outer integrators of a given level $\ell$ integrator have arrived at and is necessary to take into account to correctly reflect the numerical time of the solution on each level $\ell$.

Starting from a computed numerical solution $\fSDSys^n$ at time $t^n=nh_L$, one first takes $K_0+1$ steps of size $h_0$ with the innermost integrator,
\begin{equation} \label{eq:semidiscrete_innermost_scheme}
	\fSDSys^{0,n,k_0+1} =S_{0}(\fSDSys^{0,n,k_0}),\qquad 0 \le k_0 \le K_0, 
\end{equation}
in which $\fSDSys^{0,n,k_0}$ corresponds to the numerical solution at time $t^{0,n,k_0}$ calculated by the innermost integrator. Since all outer integrator iteration indices $k_{\ell'}$, $\ell'=1,...,L-1$ are zero in equation \eqref{eq:tpi_time}, we have ${t^{0,n,k_0} = nh_L+k_0h_0}$. The repeated action of the innermost integrator is depicted by small black arrows in the upper row of figure \ref{fig:tpi_sketch}, for which we chose $K_0=2$.

In the telescopic projective integration framework, the scheme is set up from the lowest level up to the highest level. The aim is to obtain a discrete derivative to be used on each level to eventually compute $\fSDSys^{n+1} = \fSDSys^{0,n+1,0}$ via extrapolation in time. Using the innermost integrator iterations \eqref{eq:semidiscrete_innermost_scheme}, we perform the extrapolation by a projective integrator on level 1, written as:
\begin{equation} \label{eq:TPFE_level1} 
	\fSDSys^{1,n,1} = \fSDSys^{0,n,K_{0}+1} + \left(M_{0}h_{0}\right)\frac{\fSDSys^{0,n,K_{0}+1} - \fSDSys^{0,n,K_{0}}}{h_{0}},
\end{equation}
which corresponds to the projective forward Euler (PFE) method \cite{Gear2003projective}. In equation \eqref{eq:TPFE_level1}, $\fSDSys^{1,n,1}$ represents the numerical solution at time $t^{1,n,1}$ calculated by one iteration of the first level projective integrator. Since $k_1=1$ and all its outer integrator iteration indices $k_{\ell'}$, ${\ell'=2,...,L-1}$ are still zero in equation \eqref{eq:tpi_time}, we have ${t^{1,n,1} = nh_L + h_1}$. One such step of the first level integrator is visualized by a large green arrow in the upper row of figure \ref{fig:tpi_sketch}. We can repeat this idea and construct a hierarchy of projective integrators on levels $\ell=1,...,L-1$, given by:
\begin{equation} \label{eq:TPFE_inners} 
	\fSDSys^{\ell,n,k_\ell+1} = \fSDSys^{\ell-1,n,K_{\ell-1}+1} + \left(M_{\ell-1}h_{\ell-1}\right)\frac{\fSDSys^{\ell-1,n,K_{\ell-1}+1} - \fSDSys^{\ell-1,n,K_{\ell-1}}}{h_{\ell-1}}, \qquad 0 \le k_{\ell} \le K_{\ell}. 
\end{equation}
where $\fSDSys^{\ell,n,k_\ell}$ denotes the numerical solution at time $t^{\ell,n,k_\ell}$ calculated by projective integrator on level $\ell$. According to equation \eqref{eq:tpi_time}, this time depends on the values $k_{\ell'}$, $\ell'=\ell+1,...,L-1$ of all of its outer integrators.
For each level $\ell=1,...,3$, these projective integrator steps are shown in figure \ref{fig:tpi_sketch} by long arrows. Ultimately, the outermost integrator on level $L$ computes $\fSDSys^{n+1}$ as:
\begin{equation} \label{eq:TPFE_outermost}
	\fSDSys^{n+1} = \fSDSys^{L-1,n,K_{L-1}+1} + (M_{L-1}h_{L-1})\frac{\fSDSys^{L-1,n,K_{L-1}+1} - \fSDSys^{L-1,n,K_{L-1}}}{h_{L-1}}  .
\end{equation}
Since the outermost integrator \eqref{eq:TPFE_outermost} also constitutes a PFE scheme, the telescopic method resulting from the hierarchy of projective levels \eqref{eq:TPFE_inners}-\eqref{eq:TPFE_outermost} is called telescopic projective forward Euler (TPFE), and it is the simplest instantiation of this class of integration methods. 

As shown in \cite{Lejon2016} and \cite{Melis2015}, it is straightforward to implement higher-order extensions of the outermost integrator, such as the projective Runge-Kutta methods of order 2 and 4  in the telescopic case, leading to TPRK2 and TPRK4 methods. In general, the outermost integrator in a TPRK method replaces each time derivative evaluation in a classical Runge-Kutta method, denoted by $\mathbf{k}_s$, by $K_{L-1}+1$ steps of its inner integrator on level $L-1$. Using equation \eqref{eq:TPFE_inners} with $\ell=L-1$, the first stage in a TPRK method calculates the time derivative $\mathbf{k}_1$ as:
\begin{equation} \label{eq:TPRK_stage_1}
	\mathbf{k}_1 = \dfrac{\fSDSys^{L-1,n,K_{L-1}+1} - \fSDSys^{L-1,n,K_{L-1}}}{h_{L-1}}.
\end{equation}
Any other stage $s \ge 2$ requires evaluating the time derivatives at intermediate times denoted by ${t^{n+c_s} = (n+c_s)h_L}$. Similarly to equation \eqref{eq:TPRK_stage_1}, these are calculated as:
\begin{equation} \label{eq:TPRK_stage_s}
	\mathbf{k}_s = \dfrac{\fSDSys^{L-1,n+c_s,K_{L-1}+1} - \fSDSys^{L-1,n+c_s,K_{L-1}}}{h_{L-1}}.
\end{equation}
Since the numerical solution at $t^{n+c_s}$ in equation \eqref{eq:TPRK_stage_s} is not available, we use the integrator on level $L-1$ to approximate it as follows:
\begin{equation} \label{eq:TPRK_intermediate_solutions}
	\begin{dcases}
		\fSDSys^{L-1,n+c_s,0} &= \fSDSys^{L-1,n,K_{L-1}+1} + (c_sh_L-(K_{L-1}+1)h_{L-1}) \sum_{m=1}^{s-1}\dfrac{a_{s,m}}{c_s} \mathbf{k}_s \\
		\fSDSys^{L-1,n+c_s,k_{L-1}+1} &= \fSDSys^{L-2,n+c_s,K_{L-2}+1} + \left(M_{L-2}h_{L-2}\right)\frac{\fSDSys^{L-2,n+c_s,K_{L-2}+1} - \fSDSys^{L-2,n+c_s,K_{L-2}}}{h_{L-2}},
	\end{dcases}
\end{equation}
in which the last equation in \eqref{eq:TPRK_intermediate_solutions} iterates over $0 \le k_{L-1} \le K_{L-1}$.
Ultimately, the outermost integrator of a TPRK method is written as:
\begin{equation}
	\fSDSys^{n+1} = \fSDSys^{L-1,n,K_{L-1}+1} + (M_{L-1}h_{L-1})\sum_{s=1}^{S}b_s \mathbf{k}_s.
\end{equation}
To ensure consistency, the RK matrix $\mathbf{a}=(a_{s,m})_{s,m=1}^S$, weights $\mathbf{b}=(b_s)_{s=1}^S$, and nodes ${\mathbf{c}=(c_s)_{s=1}^S}$ satisfy (see, e.g., \cite{Hairer1993}) the conditions $0\le b_s \le 1$ and $0 \le c_s \le 1,$ as well as
\begin{equation} \label{eq:RK_conditions}
	\sum_{s=1}^Sb_s=1, \qquad \sum_{m=1}^{S-1} a_{s,m} =c_s, \quad 1 \le s \le S. 
\end{equation}
(Note that these assumptions imply that $c_1=0$ using the convention that $\sum_{1}^0\cdot=0$.)

Then, the TPRK2 and TPRK4 methods are obtained by choosing their coefficients as shown in the Butcher tableaux in figure \ref{tab:butcher}. In the numerical experiments, we will specifically use the projective Runge-Kutta method of order 4.

\begin{figure}[t]
	\begin{center}
		\figname{butcher_tableaux}
		\newcommand{\butcherVerSpacing}{0.6}
\newcommand{\butcherHorSpacing}{0.8}
\newcommand{\butcherMidXstart}{4}
\newcommand{\butcherRightXstart}{9}
\begin{tikzpicture}
	\draw[](-\butcherHorSpacing/2,\butcherVerSpacing/2) -- (1.5*\butcherHorSpacing+0.15,\butcherVerSpacing/2);
	\draw[](\butcherHorSpacing/2,-\butcherVerSpacing/2) -- (\butcherHorSpacing/2,1.5*\butcherVerSpacing);
	
	\node(c) at (0,\butcherVerSpacing) {$\mathbf{c}$};
	\node(a) at (\butcherHorSpacing,\butcherVerSpacing) {$\mathbf{a}$};
	\node(bT) at (\butcherHorSpacing+0.15,0) {$\mathbf{b}^T$};
	
	\draw[] (\butcherMidXstart-\butcherHorSpacing/2,\butcherVerSpacing/2) -- (\butcherMidXstart+2.5*\butcherHorSpacing,\butcherVerSpacing/2);
	\draw[] (\butcherMidXstart+\butcherHorSpacing/2,-\butcherVerSpacing/2) --(\butcherMidXstart+\butcherHorSpacing/2,2.5*\butcherVerSpacing);
	
	\node(c1_sec) at (\butcherMidXstart,\butcherVerSpacing) {$1/2$};
	\node (c2_sec) at (\butcherMidXstart,2*\butcherVerSpacing){$0$};
	\node (a21_sec) at (\butcherMidXstart+\butcherHorSpacing,\butcherVerSpacing) {$1/2$};
	\node (b1_sec) at (\butcherMidXstart+\butcherHorSpacing,0) {$0$};
	\node (b2_sec) at (\butcherMidXstart+2*\butcherHorSpacing,0) {$1$};
	
	\draw[] (\butcherRightXstart-\butcherHorSpacing/2,\butcherVerSpacing/2) -- (\butcherRightXstart+4.5*\butcherHorSpacing,\butcherVerSpacing/2);
	\draw[] (\butcherRightXstart+\butcherHorSpacing/2,-\butcherVerSpacing/2) -- (\butcherRightXstart+\butcherHorSpacing/2,4.5*\butcherVerSpacing);

	\node (c1_fourth) at (\butcherRightXstart,\butcherVerSpacing) {$1$};	
	\node (c2_fourth) at (\butcherRightXstart,2*\butcherVerSpacing) {$1/2$};
	\node (c3_fourth) at (\butcherRightXstart,3*\butcherVerSpacing) {$1/2$};
	\node (c4_fourth) at (\butcherRightXstart,4*\butcherVerSpacing) {$0$};
	\node (a41_fourth) at (\butcherRightXstart+\butcherHorSpacing,\butcherVerSpacing) {$0$};
	\node (a42_fourth) at (\butcherRightXstart+2*\butcherHorSpacing,\butcherVerSpacing) {$0$};
	\node (a43_fourth) at (\butcherRightXstart+3*\butcherHorSpacing,\butcherVerSpacing) {$1$};
	\node (a31_fourth) at (\butcherRightXstart+\butcherHorSpacing,2*\butcherVerSpacing) {$0$};
	\node (a32_fourth) at (\butcherRightXstart+2*\butcherHorSpacing,2*\butcherVerSpacing) {$1/2$};
	\node (a21_fourth) at (\butcherRightXstart+\butcherHorSpacing,3*\butcherVerSpacing) {$1/2$};
	\node(b1_fourth) at (\butcherRightXstart+\butcherHorSpacing,0) {$1/6$};
	\node (b2_fourth) at (\butcherRightXstart+2*\butcherHorSpacing,0) {$1/3$};
	\node (b3_fourth) at (\butcherRightXstart+3*\butcherHorSpacing,0) {$1/3$};
	\node (b4_fourth) at (\butcherRightXstart+4*\butcherHorSpacing,0) {$1/6$};
\end{tikzpicture}
	\end{center}
	\vspace{-0.4cm}\caption {\label{tab:butcher} Butcher tableaux for Runge-Kutta methods. Left: general notation; middle: RK2 method (second order); right: RK4 method (fourth order).}
\end{figure}

\subsection{TPI method parameters} \label{subsec:tpi_params}
In general, the level-$L$ TPI method possesses a set of $3L+1$ parameters corresponding to all of its levels. The innermost integrator has only one parameter, its time step $h_0$. The $L$ projective integrator levels which are built around the innermost integrator each contain 3 parameters. For projective levels $\ell=1,...,L$ these are: (i) the time step $h_\ell$, (ii) the number of lower level integrator iterations $K_{\ell-1}$ that are needed to sufficiently damp the fast components at level $\ell-1$, and (iii) the extrapolation step size $M_{\ell-1}h_{\ell-1}$ over which the integrator on level $\ell-1$ is applied in the $\ell$-the level projective step. (Note that the parameter $M_{\ell-1}$ is the number of steps at level $\ell-1$ that are skipped by the projective step.) 

Notably, the innermost integrator time step $h_0$ is the only real time step of a TPI method, meaning that this is the only time step over which numerical integration is actually performed. All higher level time steps are merely the consequence of the $K$- and $M$-values of the extrapolation. Once $h_0$ is known, the projective time step on level $\ell=1,...,L$ satisfies the following relation:
\begin{equation} \label{eq:tpi_successive_timesteps} 
	h_\ell = \prod_{k=0}^{\ell-1} (M_k + K_k + 1)h_0.
\end{equation}
This can also be seen in figure \ref{fig:tpi_sketch} for a level-3 TPI method in which we chose $K=2$ constant on all levels.

\subsection{Stability of telescopic projective integration} \label{subsec:tpi_stability}
We now briefly discuss the main stability properties of the TPFE method which can be found in more detail in \cite{Gear2003telescopic}. To that end, we introduce the test equation and its corresponding innermost integrator:
\begin{equation}\label{eq:test_eq} 
	\dot{y}=\lambda y, \qquad y^{k+1}=\sigma_0(\lambda h_0)y^k, \qquad \lambda \in \C. 
\end{equation}
As in \cite{Gear2003projective}, we call $\sigma_0(\lambda h_0)$ the \emph{amplification factor} of the innermost integrator. (For instance, if the innermost integrator is the forward Euler scheme, we have $\sigma_0(\lambda h_0)=1+\lambda h_0$.) The innermost integrator is stable if $\left|\sigma_0\right|\le 1$. The question then is for which subset of these $\sigma_0$-values, which are also called $\sigma_0$-eigenvalues, the TPFE method is also stable. Considering the level-$L$ TPFE method, it can easily be seen from equations \eqref{eq:TPFE_inners}-\eqref{eq:TPFE_outermost} that it is stable if
\begin{equation} \label{eq:stab_cond} 
	|\sigma_L(\sigma_0)| = \left|\Big((M_{L-1}+1)\sigma_{L-1}(\sigma_0) - M_{L-1}\Big)\Big(\sigma_{L-1}(\sigma_0)\right)^{K_{L-1}} \Big|\le 1, 
\end{equation}
in which $\sigma_L$ denotes the outermost integrator amplification factor. Equation \eqref{eq:stab_cond} needs to hold for all eigenvalues $\sigma_0$ of the innermost integrator. 

Since we are interested in the limit $\hydroLimit$ for fixed $\dx$, we look at the limiting stability regions that arise when taking the limit $h_0 \to 0$, while keeping $h_\ell$, $\ell=1,...,L$ fixed. In this regime, it was obtained in \cite{Gear2001} that the level-$L$ TPI method contains $L+1$ (principal) regions of stability around the real axis which depend on the choice of (possibly) different $K_{\ell-1}$- and $M_{\ell-1}$-values at each projective level $\ell=1,...,L$ and can be positioned to cover the clusters of eigenvalues. Furthermore, there are a number of artefact stability regions due to the value of $K_{\ell-1}$ which can not be tuned independently and are of no importance.

The TPI method allows for the accurate integration of solution modes within its dominant (rightmost) stability region while maintaining stability for all other modes by matching its stability regions around the eigenvalue clusters of the problem's spectrum.


\section{Numerical properties} \label{sec:num_props}
In this section we describe the selection procedure of the level-$L$ TPI method parameters $h_0$, $K_{\ell-1}$ and $M_{\ell-1}$ with $\ell=1,...,L$ such that the TPI method is stable. This procedure is based on the spectrum of the innermost integrator which is derived in section \ref{subsec:spectrum_innermost}. Then, we explain the parameter selection procedure in case of a time-invariant relaxation profile $\omega(x)$ containing a number of discrete $\omega$-levels (section \ref{subsec:tpi_params_omega}) and a density-dependent relaxation time model with a time-varying spectrum (section \ref{subsec:tpi_params_rho}).

\subsection{Spectrum of the innermost integrator} \label{subsec:spectrum_innermost}
Once the spectrum of the semi-discrete system \eqref{eq:semidiscrete} is known, the spectrum of the innermost integrator can be derived by transforming its expression to the Fourier domain. For instance, when choosing the forward Euler scheme
\eqref{eq:FE_scheme} as innermost integrator, its expression in the Fourier domain is given by:
\begin{equation} \label{eq:fourier_fe}
	\fF^{k+1} = \SzeroF\;\fF^{k} = \left(\eye + h_0\tilde\Bcal \right)\fF^k,
\end{equation}
with $\SzeroF$ the Fourier transform of the forward Euler time stepper $S_0$. The matrix $\tilde{\Bcal}$ in equation \eqref{eq:fourier_fe} is either given by equation \eqref{eq:spectrum_Btilde} in case of a time-invariant collision frequency ${\nu(x,t) = \omega(x)}$ or by equation \eqref{eq:spectrum_Btilde_rho} when considering the time-varying case ${\nu(x,t) = \rho(x,t)}$.
It is clear that the amplification factors $\sigmazerovec = \set{\sigma_{0}^{(j)}}{j}$ of the forward Euler scheme, which are the eigenvalues of $\SzeroF$, and the eigenvalues $\lambdavec = \set{\lambda^{(j)}}{j}$ of the matrix $\tilde\Bcal$ are related via 
\begin{equation} \label{eq:tau_fe}
	\sigma_{0}^{(j)} = 1 + h_0\lambda^{(j)}, \qquad j = 1, ..., J.
\end{equation}
By convention, we consider the dominant eigenvalue $\lambda^{(1)}(\omegabar_1)=\lambda^{(1)}$ with $\lambda^{(1)}(\omegabar_1)$ given in theorem \ref{thm_spec}.

\subsection{Time-invariant collision frequency} \label{subsec:tpi_params_omega}
In this setting, the starting point is a given eigenvalue spectrum containing $\Lw+1$ eigenvalue clusters with $\Lw$ fast and 1 slow clusters, which are all located in the left half plane of the complex $\lambda$-plane. In this case, the $\Lw$ fast clusters, which may or may not be clearly separated, arise from $\Lw$ different $\omegabar_{\lw}$-values, $\lw=1, ..., \Lw$ in the piecewise constant relaxation function $\omega(x)$. We label the eigenvalue cluster centers in the $\lambda$-plane such that ${\lambda_0 \le \lambda_1 \le \cdots \le \lambda_{\Lw} \le 0}$. This means that $\lambda_0$ corresponds to the center of the fastest eigenvalue cluster whereas $\lambda_{\Lw}$ represents the cluster center with dominant (slow) eigenvalues. We only assume a clear spectral gap between the fastest and slow cluster, i.e., $\lambda_0 \ll \lambda_{\Lw}$, such that it is useful to implement a projective method. In what follows, we first elucidate the underlying idea of choosing the TPI method parameters. Afterwards, we provide a detailed description.

For all levels $\ell$ from 0 to $L-1$, which contain an integrator that serves an inner integrator role, the purpose of the integrator on level $\ell$ is to bring the $\ell^{th}$ fast eigenvalue cluster to 0 (i.e., the integrator on level $\ell$ should damp all eigenvalues in the $\ell^{th}$ fast cluster). Notice that we use $\{\ell,L\}$ for levels of the TPI method and $\{\lw,\Lw\}$ for values in the relaxation function $\omega(x)$ which are not necessarily the same (see later). The $\ell$ faster-than-the-current eigenvalue clusters, which were already around 0 due to the application of the lower level integrators, remain around 0 and are suppressed even more by the integrator on level $\ell$. The $L-\ell$ remaining eigenvalue clusters to the right of the current ($\ell^{th}$) eigenvalue cluster will shift somewhat more to the left (i.e., towards $0$) since the integrator on level $\ell$ is also (slightly) damping these eigenvalues. For stability reasons, we require that the $\ell$ clusters that were already around 0 lie in the stability region of the integrator on level $\ell+1$. This can be achieved by carefully selecting the value of $K_{\ell-1}$ on each level.
It is clear that in the above reasoning we only desire a stable numerical integration of the fast modes. However, the dominant modes, which are the solution components of practical interest, need to be integrated both in a stable and accurate way. 

Below, we detail the general level-$L$ TPI method construction, for which we always assume the forward Euler method on levels $\ell=0,...,L-1$. The outermost ($L^{th}$ level) integrator can be any stable explicit method, depending on the required accuracy. We distinguish between two cases. First, in section \ref{subsubsec:spectral_gaps}, we assume that the eigenvalue clusters are clearly separated. Afterwards, in section \ref{subsubsec:no_spectral_gaps}, we comment on the situation when this assumption is not satisfied. In both cases, we discuss the construction of TPI methods and illustrate with numerical results.

\subsubsection{Spectrum with spectral gaps} \label{subsubsec:spectral_gaps}
\begin{figure}[t]
	\begin{center}
		\figname{tpfe_spectrum_omega_2levels}
		\input{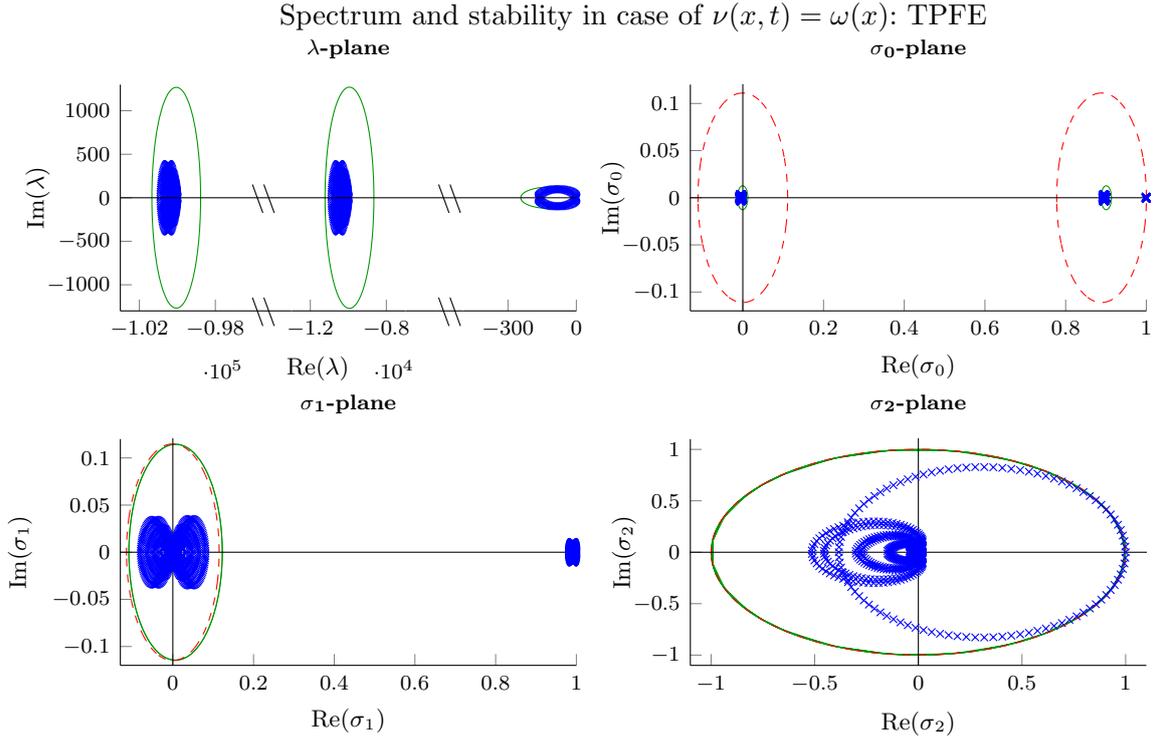}
	\end{center}
	\vspace{-0.4cm}\caption{\label{fig:spectrum_tpfe_omega} Stability analysis for a level-2 TPFE method in case of a time-invariant spectrum with clear spectral gaps. Blue crosses in each plane correspond to eigenvalues in that plane. Green regions represent the stability regions of the full level-2 TPFE method with respect to every particular plane. The red dashed regions are the stability regions of the next level projective integrator in a certain plane. }
\end{figure}

As we proceed through the detailed construction procedure below we use the following numerical experiment. We fix $\epsi=10^{-5}$ and discretize velocity space using $J=10$ velocity components obtained as the nodes of Gauss-Hermite quadrature for integration with respect to the measure given in \eqref{eq:V_measure}. We consider $x\in[0, 1]$, together with periodic boundary conditions and propose a piecewise constant relaxation profile $\omega(x)$ containing 2 well separated $\omega$-values: $\{1,0.1\}$ such that $\omega(x) = 1\;(x \le 0.5) $ and $\omega(x) = 0.1\;(x > 0.5)$. We use the upwind scheme of order $1$ with grid spacing $\dx=0.01$ as spatial discretization technique. The resulting spectrum of system \eqref{eq:semidiscrete} in the $\lambda$-plane is shown by blue crosses in the top left plot of figure \ref{fig:spectrum_tpfe_omega}. As can be seen, there are 2 fast clusters centered around the positions $-1/\epsi$ and $-0.1/\epsi$ as formalized in equation \eqref{eq:spectrum_Btilde}. Since there are 2 clearly separated fast clusters, we construct a level-2 TPI method. In general, given that we consider clearly separated clusters in this section, we have $L=\Lw$.

\paragraph{Innermost integrator.}
The innermost integrator of the TPI method corresponds to a space-time discretization of equation \eqref{eq:kin_eq_bgk_1d}, in which we choose the simple first-order explicit forward Euler time discretization with time step $h_0$. As explained in the introduction of this section, we fix $h_0$ such that the fastest eigenvalue cluster with center $\lambda_0$ is moved to 0 in the complex $\sigma_0$-plane. Using equation \eqref{eq:tau_fe} and the result in equation \eqref{eq:spectrum_Btilde}, we obtain:
\begin{equation} \label{eq:fe_h0} 
	\sigma_0(\lambda_0) = 0 \Rightarrow h_0 = \frac{1}{|\lambda_0|} = \frac{\epsi}{\omegabar_1},
\end{equation}
This choice of the time step $h_0$ defines a linear mapping of all eigenvalues $\lambda$ to (zeroth level) eigenvalues $\sigma_0$, which can all be found in the $\sigma_0$-plane within the interval $[-\eta,1]$ with $\eta \in \R^{+}$ close to zero. 

In the numerical example, given that $\omegabar_1=1$, equation \eqref{eq:fe_h0} gives rise to $h_0=\epsi$. The corresponding $\sigma_0$-eigenvalues are visualized by blue crosses in the top right plot of figure \ref{fig:spectrum_tpfe_omega}.

\paragraph{Integrator on level $\ell$.}
For $\ell=1,...,L-1$ the integrator on level $\ell$ is the projective forward Euler (PFE) scheme, which extrapolates the solution of its inner integrator (i.e. the integrator on level $\ell-1$) over a distance $M_{\ell-1}h_{\ell-1}$. Its amplification factor $\sigma_\ell$ in terms of its inner integrator amplification factor $\sigma_{\ell-1}$ is given by:
\begin{equation} \label{eq:pfel_amp_fact} 
	\sigma_{\ell}(\sigma_{\ell-1}) = \Big((M_{\ell-1}+1)\sigma_{\ell-1} - M_{\ell-1}\Big)(\sigma_{\ell-1})^{K_{\ell-1}}. 
\end{equation}
After applying the integrator on level $\ell-1$, there will be $L+2-\ell$ eigenvalue clusters remaining in the $\sigma_{\ell-1}$-plane of which there are $L+1-\ell$ fast and 1 slow cluster, and $\ell$ clusters will already have been moved to 0 by the lower level integrators turning these into one big cluster around 0. 
We then demand that $M_{\ell-1}$ is chosen such that the fastest eigenvalue cluster that is not yet around 0 in the $\sigma_{\ell-1}$-plane will be moved to 0 in the $\sigma_\ell$-plane. This cluster is denoted by $\sigma_{\ell-1,\ell}$, in which the notation $\sigma_{m,n}=\sigma_m(\lambda_n)$ represents the transformation of eigenvalue cluster $\lambda_n$ to the $\sigma_m$-plane. Using expression \eqref{eq:pfel_amp_fact} we find:
\begin{equation} \label{eq:pfel_Ml-1} 
	M_{\ell-1} \mbox{ such that } \sigma_\ell(\sigma_{\ell-1,\ell}) = 0 \Rightarrow M_{\ell-1} = \frac{\sigma_{\ell-1,\ell}}{1 - \sigma_{\ell-1,\ell}}. 
\end{equation}
Note that this choice of $M_{\ell-1}$ makes the integrator on level $\ell$ maximally damping at the center of eigenvalue cluster $\sigma_{\ell-1,\ell}$ (i.e. it moves the center of cluster $\sigma_{\ell-1,\ell}$ exactly to 0). We point out that the extrapolation step sizes $M_{\ell-1}$ in equation \eqref{eq:pfel_Ml-1} only depend on the ratio of two consecutive eigenvalue clusters in the $\lambda$-plane, and are thus independent of $\epsi$. This can be shown by working out equation \eqref{eq:pfel_Ml-1} employing the same Taylor expansion as in equation \eqref{eq:pfeL_sigl_fd} and using the result in equation \eqref{eq:spectrum_Btilde} yielding:
\begin{align} \label{eq:M_no_epsi} 
	M_{\ell-1} = \frac{\sigma_{\ell-1,\ell}}{1 - \sigma_{\ell-1,\ell}} &\approx \frac{1 + h_{\ell-1}\lambda_{\ell}}{-h_{\ell-1}\lambda_\ell} \notag \\
	& = \frac{\epsi}{h_{\ell-1}\omegabar_{\ell+1}} - 1 \notag \\
	& \approx \frac{\omegabar_\ell}{\omegabar_{\ell+1}} - 1,
\end{align}
where we used in the last step that the time step $h_{\ell-1} = O(\epsi/\omegabar_\ell)$ commensurate with the time scale of the previous eigenvalue cluster. 

Once the value of $M_{\ell-1}$ is known, we fix the value of $K_{\ell-1}$ by demanding that all fast eigenvalue clusters that are already around 0 in the $\sigma_{\ell-1}$-plane fall into the stability region of the projective integrator on level $\ell$ around 0, which is given by $\disk(0,(1/M_{\ell-1})^{1/K_{\ell-1}})$ \cite{Gear2003projective} where $\disk(c,r)$ denotes the disk with center $(c,0)$ and radius $r$. We thus obtain:
\begin{equation} \label{eq:pfel_Kl-1} 
	\hat\sigma_{\ell-1} \le \left(\frac{1}{M_{\ell-1}}\right)^{1/K_{\ell-1}} \Rightarrow K_{\ell-1} = \left\lceil\frac{\log(1/M_{\ell-1})}{\log(\hat\sigma_{\ell-1})}\right\rceil,
\end{equation}
in which $\hat\sigma_{\ell-1} = \max_{k\in\{0,1,...,\ell-1\}} |\sigma_{\ell-1,k}|$ denotes the cluster that is farthest away from the origin. Since the values $M_{\ell-1}$ are independent of $\epsi$, we find that the values $K_{\ell-1}$ in equation \eqref{eq:pfel_Kl-1} are also independent of $\epsi$. 

In the numerical experiment, using equations \eqref{eq:pfel_Ml-1} and \eqref{eq:pfel_Kl-1} we construct the first level PFE method with parameters $M_0=9$ and $K_0=1$. The (asymptotic, $M_0 \to \infty$) stability regions corresponding to this projective level integrator are indicated by red dashed circles in the top right plot of figure \ref{fig:spectrum_tpfe_omega}. The eigenvalues transformed to the $\sigma_1$-plane are shown by blue crosses in the bottom left plot of figure \ref{fig:spectrum_tpfe_omega}.

\paragraph{Outermost integrator.}
Finally, the integrator on level $L$ or outermost integrator is designed by ensuring that all eigenvalues belonging to the dominant eigenvalue cluster  $\sigma_{L-1,L}$ in the $\sigma_{L-1}$-plane fall in the dominant stability region of the outermost integrator, which is given by the region $\disk(1-1/M_{L-1},1/M_{L-1})$ \cite{Gear2003projective}. This leads to the following inequality:
\begin{equation} \label{eq:pfeL_condition} 
	\sqrt{\left(\Re(\sigma_{L-1,L}) -  \left(1-\frac{1}{M_{L-1}}\right)\right)^2 + \Big(\Im(\sigma_{L-1,L})\Big)^2} \le \frac{1}{M_{L-1}}.
\end{equation}
To calculate the value of $M_{L-1}$ from \eqref{eq:pfeL_condition}, we examine how the expression of the dominant eigenvalues $\sigma_{0,L}$ in the $\sigma_0$-plane, which will remain the dominant eigenvalues in all other $\sigma_\ell$-planes, $\ell=1,...,L$, transforms under application of the different level projective integrators.
Notice that this expression is already known in the $\sigma_0$-plane using equations \eqref{eq:spectrum_Btilde} and \eqref{eq:tau_fe}, and is of the following form:
\begin{equation} \label{eq:pfeL_sig0_fd} 
	\sigma_{0,L} = \sigma_0(\lambda_L) = 1 + h_0\lambda_L. 
\end{equation}
By plugging the expression for $\sigma_{0,L}$ in \eqref{eq:pfeL_sig0_fd} into equation \eqref{eq:pfel_amp_fact} we find that the first level PFE integrator scheme transforms these dominant $\sigma_0$-eigenvalues into dominant $\sigma_1$-eigenvalues:
\begin{equation} \label{eq:pfeL_sig1_fd} 
	\sigma_{1,L} \approx 1 + h_{1}\lambda_L.
\end{equation}
The expression for $\sigma_{1,L}$ is obtained by performing a Taylor series expansion for $h_0 \to 0$ and using the time step relation given in equation \eqref{eq:tpi_successive_timesteps}.
Repeating this line of thought for the next level integrators, we find that for $\ell=1,...,L-1$ the dominant eigenvalues transformed by the integrator on level $\ell$ is written as follows:
\begin{equation} \label{eq:pfeL_sigl_fd} 
	\sigma_{\ell,L} \approx 1 + h_{\ell}\lambda_L.
\end{equation}
Using equation \eqref{eq:pfeL_sigl_fd}, we can now use the condition given in \eqref{eq:pfeL_condition} to find the value of $M_{L-1}$. To this end, we will turn the inequality in \eqref{eq:pfeL_condition} into an equality. Using equations \eqref{eq:dom_eig_real}-\eqref{eq:dom_eig_imag}, the corresponding (maximum allowed) value of $M_{L-1}$ is given by
\begin{equation} \label{eq:pfeL_M} 
	M_{L-1} = \min_\zeta \left(\frac{-2\averageV{\alphavec}}{h_{L-1}\Big(\averageV{\alphavec}^2 + \averageV{\betavec\vGrid}^2\Big)}\right).
\end{equation}
Once the value of $M_{L-1}$ is known, we use equation \eqref{eq:pfel_Kl-1} to determine the corresponding value of $K_{L-1}$.

From equation \eqref{eq:pfeL_M}, we observe that the extrapolation step size $M_{L-1}$ of the outermost integrator, which bridges the gap between the last fast cluster and the dominant slow cluster, is inversely proportional to $h_{L-1}$. Given that the latter depends on $\epsi$, see equation \eqref{eq:tpi_successive_timesteps}, we conclude that $M_{L-1}$ is inversely proportional to $\epsi$ as desired.

In the numerical example, using equations \eqref{eq:pfeL_M} and \eqref{eq:pfel_Kl-1}, the parameters of the outermost PFE integrator are given by $M_1 = 75.82$ and $K_1 = 2$. The CFL number for this choice of parameters is $0.87$.
The (asymptotic, $M_1 \to \infty$) stability regions of the outermost projective level are visualized by red dashed circles in the bottom left plot of figure \ref{fig:spectrum_tpfe_omega}. The $\sigma_2$-eigenvalues are depicted by blue crosses in the bottom right plot of figure \ref{fig:spectrum_tpfe_omega}.  We also plotted the stability domain of the full level-2 TPFE method in each complex plane using solid green lines. In particular, notice that in the $\sigma_2$-plane this stability domain coincides with the unit disc and all $\sigma_2$-eigenvalues lie within this region.

\subsubsection{Spectrum without spectral gaps} \label{subsubsec:no_spectral_gaps}
The given piecewise constant relaxation profile $\omega(x)$ consisting of $\Lw$ different values $\omegabar_{\lw}$, ${\lw=1, ..., \Lw}$ gives rise to $\Lw$ eigenvalue clusters. However, in general, a clear spectral gap between two or more consecutive clusters does not necessarily exist. In that case, it makes more sense to combine two or more such consecutive clusters into one big cluster. Since each of these resulting big clusters introduces a projective level we relabel them by $\ell=1,...,L$ with $L < \Lw$. The criterion used here to decide upon creating a big cluster is to require a minimum value $M_{\min}$ of $M_{\ell-1}$ in equation \eqref{eq:pfel_Ml-1} on each level which serves as a measure for spectral separation of clusters. Then, if equation \eqref{eq:pfel_Ml-1} yields $M_{\ell-1} < M_{\min}$, this implies that the current and next cluster are too close to each other to be considered as two distinct clusters. Consequently, the selection procedure skips the next cluster and moves on to the following eigenvalue cluster instead. In the numerical experiments, we choose $M_{\min} = 3$.

\begin{figure}[t]
	\begin{center}
		\figname{tpfe_spectrum_omega_2levels_nogap}
		\input{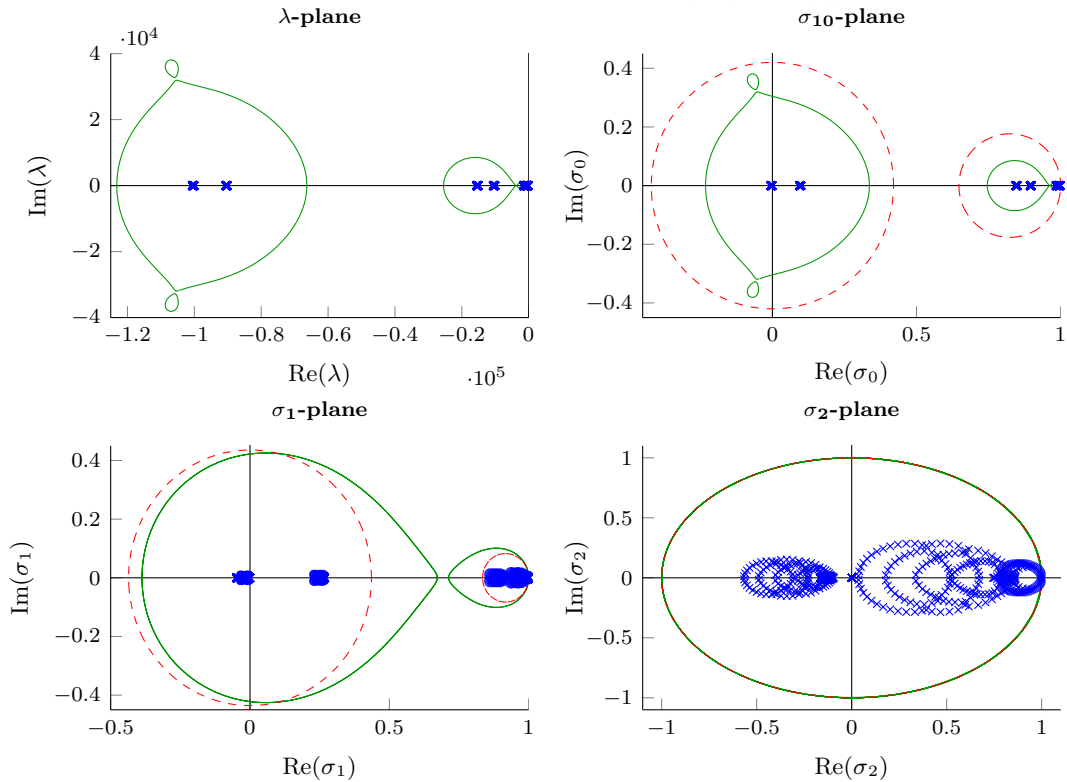}
	\end{center}
	\vspace{-0.4cm}\caption{\label{fig:spectrum_tpfe_omega_nogap} Same setup as in figure \ref{fig:spectrum_tpfe_omega}. However, not all clusters are clearly separated now. }
\end{figure}

We illustrate the construction process when there is no clear gap between every cluster. We use the same parameters as in the previous experiment. We propose a relaxation profile containing 6 $\omega$-values $\{1, 0.9, 0.15, 0.1, 0.01, 0.001\}$ which introduces 6 fast clusters in the spectrum positioned at $-\omegabar_{\lw}/\epsi$, $\lw=1,...,6$. The eigenvalues are visualized by blue crosses in the top left plot of figure \ref{fig:spectrum_tpfe_omega_nogap}. As explained above, we probe for distinct fast clusters by requiring a minimal value $M_{\min}$ of $M$ on each level. When putting $M_{\min} = 3$ we only retrieve 2 true fast clusters. Moreover, the algorithm detects that the clusters corresponding to the last two $\omega$-values should be understood as extra slow clusters since they appear rather close to the true slow cluster around 0. Therefore, we have $L=2 < 6=\Lw$. Running through the TPI construction procedure, we now obtain a stable level-2 TPFE method with parameters $M = \{5.67, 12.09\}$, $K = \{2,3\}$ and corresponding CFL number 0.14. The stability regions of the level-2 TPFE method and eigenvalues are shown in each plane by solid green lines and blue crosses, respectively. We observe that the first two stability regions each match two consecutive fast clusters in the $\lambda$- and $\sigma_0$-planes. Furthermore, in these planes, we notice two small artifact stability regions close to the leftmost stability region which are not used, see section \ref{subsec:tpi_stability}. There is also a very small stability region around the slow clusters which is hard to discern in the $\lambda$- and $\sigma_0$-planes due to the scaling used.

\begin{figure}[t]
	\begin{center}
		\figname{tprk4_spectrum_omega_4levels}
		\input{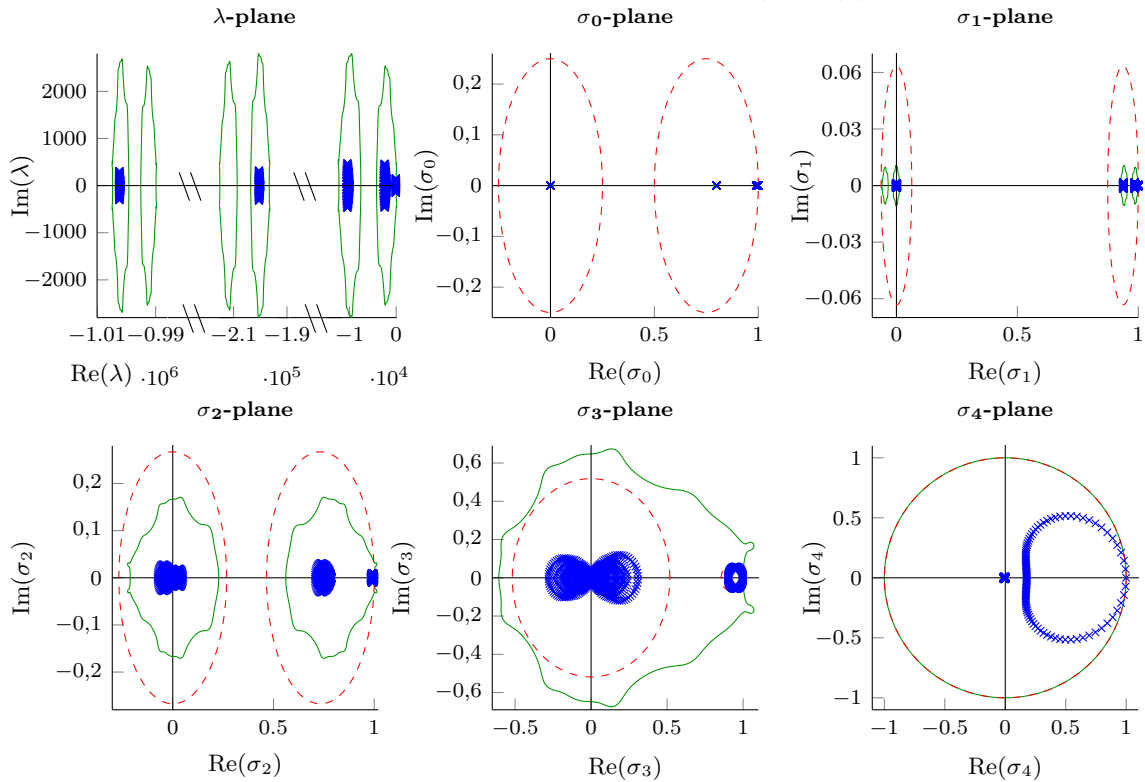}
	\end{center}
	\vspace{-0.4cm}\caption{\label{fig:spectrum_tprk4_omega} Stability analysis for a level-4 TPRK4 method in case of spectrum with clear spectral gaps. The same plotting style as in figure \ref{fig:spectrum_tpfe_omega} is used. }
\end{figure}

As a last result, we consider $\epsi=10^{-6}$ and choose 4 well separated $\omega$-levels $\{1, 0.2, 0.01, 0.002\}$. The resulting spectrum in the $\lambda$-plane, consisting of 4 fast eigenvalue clusters with clear spectral gaps and 1 slow cluster, is shown by blue crosses in the top left plot of figure \ref{fig:spectrum_tprk4_omega}. Now, we select the PFE method on each level and choose the PRK4 method as outermost integrator. As explained in \cite{Melis2015}, suitable parameters for PFE will also be suitable for higher-order PRK methods. Therefore, we determine the parameters of the outermost PRK4 integrator based on those of the PFE method. Since $L = \Lw = 4$, we construct a level-4 TPRK4 method. The parameters are determined as $M = \{4.00, 15.81, 3.74, 13.88\}$, $K = \{1, 1, 1, 4\}$ and the CFL number is 1.16. The results can be seen in figure \ref{fig:spectrum_tprk4_omega}. We observe that the constructed level-4 TPRK4 method has 5 stability regions that match the 4 fast and 1 dominant cluster.

We conclude that, in case of a time-invariant spectrum, the TPI method cost is completely independent of $\epsi$ and can be bounded from above by the number of $\omega$-levels. To illustrate this claim, consider a spectrum in the $\lambda$-plane for fixed $\epsi$ based on $\Lw$ $\omega$-levels in which some fast cluster centers $\omegabar_{\lw}/\epsi$ lie close to the dominant cluster (i.e., $\omegabar_{\lw} = O(\epsi)$ for some $\lw$). Such fast clusters are regarded as extra slow clusters and are required to lie in the dominant stability region of the TPI method. Consequently, these extra slow clusters do not introduce additional projective levels and we have $L < \Lw$. When lowering the value of $\epsi$ using the same $\omega$-profile, these extra slow clusters shift to the left one by one in the $\lambda$-plane, thus becoming true fast clusters. Each such newly obtained fast cluster requires an additional projective level to generate an additional stability region around this cluster. As soon as all  extra slow clusters became fast clusters, the number of projective levels remains constant for $\epsi \to 0$ yielding $L = \Lw$. Therefore, the cost initially behaves as $\log(1/\epsi)$ and is bounded from above by the number of $\omega$-levels which is constant and independent of $\epsi$.

\subsection{Time-varying collision frequency} \label{subsec:tpi_params_rho}
In case of a time-varying collision frequency $\collFreq(x,t) = \rho(x,t)$ the spectrum varies in principle continuously on the negative real axis. Even though the collision frequency is discretized in space we are required to take into account a continuous range of eigenvalues due to the time dependence of the collision frequency and the resulting spectrum, see equation \eqref{eq:spectrum_Btilde_rho}. In what follows, we require that the proposed numerical methods obey a maximum principle, meaning that for all discrete times $t^n$ we have:
\begin{equation} \label{eq:maximum_principle}
	0 \le \min_{1 \le i \le I}\rho^n_i, \qquad \max_{1 \le i \le I}\rho^n_i \le \max_{x}\rho(x,0),
\end{equation}
where $\rho^n_i$ is the numerical approximation of the particle density at time $t^n$ on grid point $x_i$. This maximum principle guarantees that (i) eigenvalues never cross from the left to the right side of the complex plane leading to an unstable system (left inequality), and (ii) time integration does not generate eigenvalues that become more negative than the initial fastest eigenvalue cluster (right inequality).

As set out in the introduction (section \ref{sec:intro}), the TPI method parameters can be designed such that its stability region does not split up and covers a continuous range of eigenvalues with a bigger speedup than classical projective integration.

\begin{figure}[t]
	\begin{center}
		\figname{tpfe_spectrum_rho_2levels}
%
%
\definecolor{mycolor1}{rgb}{0.00000,0.58824,0.00000}%
\newcommand{\figWidth}{5.5cm}
\newcommand{\figHeight}{3.5cm}
\newcommand{\figSpacingRight}{1.5cm} 
\newcommand{\figSpacingTop}{1.7cm} 
\begin{tikzpicture}

\begin{axis}[%
width=0,
height=0,
scale only axis,
clip=false,
xmin=0,
xmax=1,
ymin=0,
ymax=1,
hide axis
]
\node[align=center, text=black] at (\figWidth+\figSpacingRight/2,2*\figHeight+\figSpacingTop+0.9cm) {Spectrum and stability in case of $\nu(x,t)=\rho(x,t)$: TPFE};
\end{axis}

\begin{axis}[%
width=\figWidth,
height=\figHeight,
at={(0,\figHeight+\figSpacingTop)},
scale only axis,
xmin=-1.58e5,
xmax=10,
xlabel={$\text{Re(}\lambda\text{)}$},
ymin=-6.4e4,
ymax=6.4e4,
ylabel={$\text{Im(}\lambda\text{)}$},
axis background/.style={fill=white},
title style={font=\scriptsize\bfseries},
title={$\lambda\text{-plane}$},
axis x line*=bottom,
axis y line*=left
]
\addplot [color=blue,only marks,mark=x,mark options={solid},forget plot]
  table[]{tikz/data/spectrum_tpfe_rho_2levels-1.tsv};
\addplot [color=black,solid,forget plot]
  table[]{tikz/data/spectrum_tpfe_rho_2levels-2.tsv};
\addplot [color=black,solid,forget plot]
  table[]{tikz/data/spectrum_tpfe_rho_2levels-3.tsv};
\addplot [color=mycolor1,solid,forget plot,each nth point=10]
  table[]{tikz/data/spectrum_tpfe_rho_2levels-4.tsv};
\addplot [color=mycolor1,solid,forget plot,each nth point=10]
  table[]{tikz/data/spectrum_tpfe_rho_2levels-5.tsv};
\addplot [color=mycolor1,solid,forget plot,each nth point=10]
  table[]{tikz/data/spectrum_tpfe_rho_2levels-6.tsv};
\addplot [color=mycolor1,solid,forget plot,each nth point=10]
  table[]{tikz/data/spectrum_tpfe_rho_2levels-7.tsv};
\addplot [color=black,solid,forget plot]
  table[]{tikz/data/spectrum_tpfe_rho_2levels-8.tsv};
\addplot [color=black,solid,forget plot]
  table[]{tikz/data/spectrum_tpfe_rho_2levels-9.tsv};
\end{axis}

\begin{axis}[%
width=\figWidth,
height=\figHeight,
at={(\figWidth+\figSpacingRight,\figHeight+\figSpacingTop)},
scale only axis,
xmin=-0.58,
xmax=1,
xlabel={$\text{Re(}\sigma{}_\text{0}\text{)}$},
ymin=-0.63,
ymax=0.63,
ylabel={$\text{Im(}\sigma{}_\text{0}\text{)}$},
ylabel style={yshift=-0.3cm},
axis background/.style={fill=white},
title style={font=\scriptsize\bfseries},
title={$\sigma{}_\text{0}\text{-plane}$},
axis x line*=bottom,
axis y line*=left
]
\addplot [color=blue,only marks,mark=x,mark options={solid},forget plot]
  table[]{tikz/data/spectrum_tpfe_rho_2levels-10.tsv};
\addplot [color=black,solid,forget plot]
  table[]{tikz/data/spectrum_tpfe_rho_2levels-11.tsv};
\addplot [color=black,solid,forget plot]
  table[]{tikz/data/spectrum_tpfe_rho_2levels-12.tsv};
\addplot [color=mycolor1,solid,forget plot,each nth point=10]
  table[]{tikz/data/spectrum_tpfe_rho_2levels-13.tsv};
\addplot [color=mycolor1,solid,forget plot,each nth point=10]
  table[]{tikz/data/spectrum_tpfe_rho_2levels-14.tsv};
\addplot [color=mycolor1,solid,forget plot,each nth point=10]
  table[]{tikz/data/spectrum_tpfe_rho_2levels-15.tsv};
\addplot [color=mycolor1,solid,forget plot,each nth point=10]
  table[]{tikz/data/spectrum_tpfe_rho_2levels-16.tsv};
\addplot [color=black,solid,forget plot]
  table[]{tikz/data/spectrum_tpfe_rho_2levels-17.tsv};
\addplot [color=black,solid,forget plot]
  table[]{tikz/data/spectrum_tpfe_rho_2levels-18.tsv};
\end{axis}

\begin{axis}[%
width=\figWidth,
height=\figHeight,
at={(0,0)},
scale only axis,
xmin=-0.63,
xmax=1,
xlabel={$\text{Re(}\sigma{}_\text{1}\text{)}$},
ymin=-0.67,
ymax=0.67,
ylabel={$\text{Im(}\sigma{}_\text{1}\text{)}$},
axis background/.style={fill=white},
title style={font=\scriptsize\bfseries},
title={$\sigma{}_\text{1}\text{-plane}$},
axis x line*=bottom,
axis y line*=left
]
\addplot [color=blue,only marks,mark=x,mark options={solid},forget plot]
  table[]{tikz/data/spectrum_tpfe_rho_2levels-19.tsv};
\addplot [color=black,solid,forget plot]
  table[]{tikz/data/spectrum_tpfe_rho_2levels-20.tsv};
\addplot [color=black,solid,forget plot]
  table[]{tikz/data/spectrum_tpfe_rho_2levels-21.tsv};
\addplot [color=black,solid,forget plot]
  table[]{tikz/data/spectrum_tpfe_rho_2levels-22.tsv};
\addplot [color=black,solid,forget plot]
  table[]{tikz/data/spectrum_tpfe_rho_2levels-23.tsv};
\addplot [color=mycolor1,solid,forget plot]
  table[]{tikz/data/spectrum_tpfe_rho_2levels-24.tsv};
\addplot [color=mycolor1,solid,forget plot]
  table[]{tikz/data/spectrum_tpfe_rho_2levels-25.tsv};
\end{axis}

\begin{axis}[%
width=\figWidth,
height=\figHeight,
at={(\figWidth+\figSpacingRight,0)},
scale only axis,
xmin=-1.1,
xmax=1.1,
xlabel={$\text{Re(}\sigma{}_\text{2}\text{)}$},
ymin=-1.1,
ymax=1.1,
ylabel={$\text{Im(}\sigma{}_\text{2}\text{)}$},
ylabel style={yshift=-0.3cm},
axis background/.style={fill=white},
title style={font=\scriptsize\bfseries},
title={$\sigma{}_\text{2}\text{-plane}$},
axis x line*=bottom,
axis y line*=left
]
\addplot [color=blue,only marks,mark=x,mark options={solid},forget plot]
  table[]{tikz/data/spectrum_tpfe_rho_2levels-26.tsv};
\addplot [color=black,solid,forget plot]
  table[]{tikz/data/spectrum_tpfe_rho_2levels-27.tsv};
\addplot [color=black,solid,forget plot]
  table[]{tikz/data/spectrum_tpfe_rho_2levels-28.tsv};
\addplot [color=mycolor1,solid,forget plot]
  table[]{tikz/data/spectrum_tpfe_rho_2levels-31.tsv};
\end{axis}
\end{tikzpicture}%
	\end{center}
	\vspace{-0.4cm}\caption{\label{fig:spectrum_tpfe_rho} Stability analysis for a level-2 TPFE method in case of a collision frequency $\collFreq(x,t)=\rho(x,t)$. Blue crosses in each plane correspond to eigenvalues in that plane. Green regions represent the stability regions of the full level-2 TPFE method with respect to the plane under study. Clearly, the eigenvalues are distributed along the negative axis. Consequently, the stability region is not allowed to split into multiple disks. }
\end{figure}
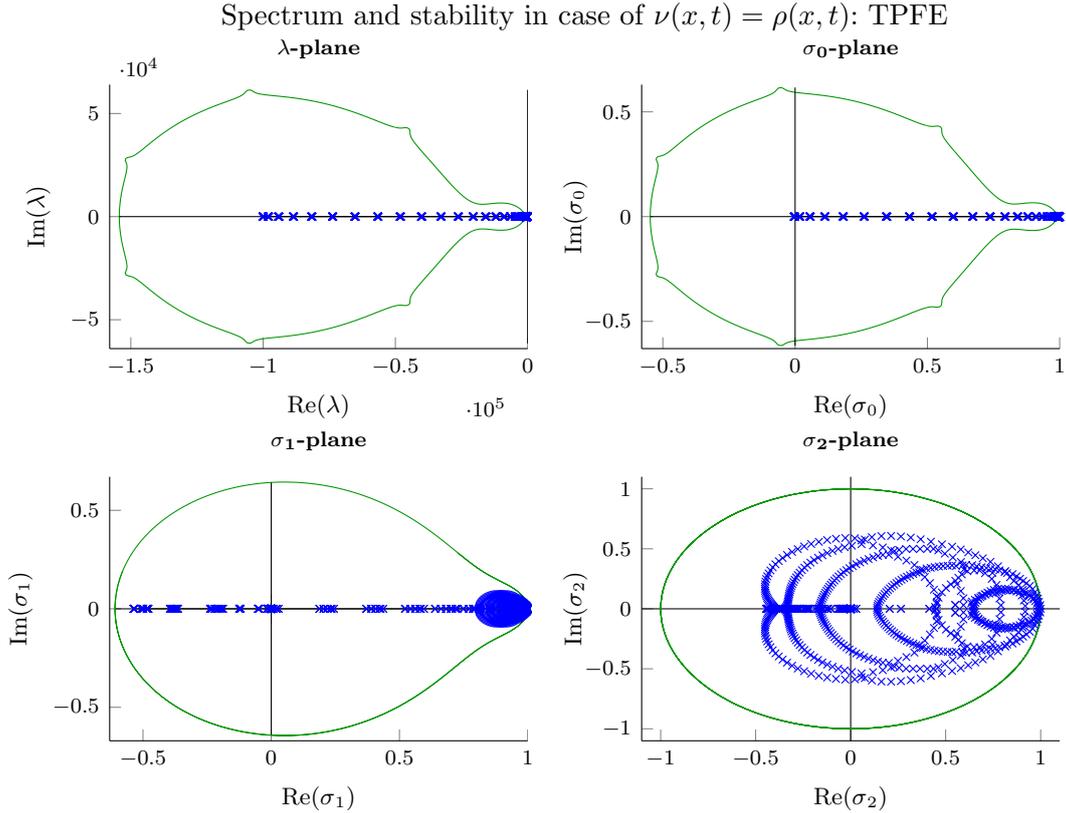

As in section \ref{subsec:tpi_params_omega}, we numerically demonstrate the method construction steps outlined below using the same setup as before. However, since the density $\rho(x,t)$ changes in time, we only visualize the spectrum corresponding to the initial condition, for which we choose the following continuous Gaussian function:
\begin{equation} \label{eq:spectrum_rho_init_rho}
	\rho(x,0) = \exp(-100(x-0.5)^2).
\end{equation}
The eigenvalues in the $\lambda$-plane are shown by blue crosses in the top left plot of figure \ref{fig:spectrum_tpfe_rho}.

Given that we use the first-order upwind scheme which satisfies both conditions in \eqref{eq:maximum_principle} we are assured that when applying the TPI construction procedure below based on the initial density in equation \eqref{eq:spectrum_rho_init_rho} we obtain a stable TPI method that remains stable for all times.

\paragraph{Innermost integrator.}
We again consider the simple first-order explicit forward Euler scheme with time step $h_0$ as innermost integrator. The requirement for choosing $h_0$ remains the same as before: it is selected such that the fastest eigenvalue in the $\lambda$-plane is moved to 0 in the $\sigma_0$-plane. In case of a continuously varying spectrum, this is done by exploiting the maximum principle \eqref{eq:maximum_principle} which guarantees that no faster cluster than the fastest cluster resulting from the initial particle density can appear during time integration. Using equations \eqref{eq:spectrum_Btilde_rho} and \eqref{eq:tau_fe}, this yields:
\begin{equation} \label{eq:fe_h0_rho} 
	h_0 = \frac{\epsi}{\max\limits_{1 \le i \le I}\rho(x_i,0)}.
\end{equation}
Using equation \eqref{eq:fe_h0_rho} we calculate the innermost integrator time step $h_0$ in the numerical experiment as $h_0=\epsi$. The eigenvalues in the $\sigma_0$-plane can be seen by blue crosses in the top right plot of figure \ref{fig:spectrum_tpfe_rho}.

As numerical diffusion may appear, this fixed choice of $h_0$ in equation \eqref{eq:fe_h0_rho} may become too restrictive for larger times $t^n$ but nevertheless is needed to integrate the first steps in a stable way. A possible extension would be to adaptively select $h_0$ depending on the maximal particle density obtained with the numerical scheme.

\paragraph{Outer integrators.}
Wrapped around the innermost integrator we construct $L$ outer integrators. Here, we always choose the projective forward Euler (PFE) scheme as outer integrator except for the outermost integrator which can be any explicit integrator as discussed before. Since we need to take care of a continuous range of eigenvalues, we calculate the method parameters $M_{\ell-1}$ and $K_{\ell-1}$ such that each outer integrator on level $\ell=1,...,L$ is $[0,1]$-stable in the $\sigma_{\ell-1}$-plane, meaning that its stability region does not split up into two disks but instead always covers the interval $[0,1]$ in that plane.

First, the required number $L$ of outer integrators is obtained by expressing that the outermost time step $h_L$ defined in equation \eqref{eq:tpi_successive_timesteps} is limited by the expected CFL stability constraint $h_L = C\dx$ for the limiting equation \eqref{eq:lin_adv}. Assuming constant values $M$ and $K$ on each level, we obtain:
\begin{equation}
	(M + K + 1)^Lh_0 \le C\dx.
\end{equation}
Consequently, the required number of projective levels is calculated as \cite{Gear2003telescopic}:
\begin{equation} \label{eq:tpi_zero_one_stable_L}
	L \approx \frac{\log(C\dx) + \log(1/h_0)}{\log(M + K + 1)}.
\end{equation}
Next, we fix the value of $K$ which is considered to be the same on each level. From the chosen value of $K$ we calculate the maximal value of $M$ needed for a $[0,1]$-stable outer integrator on each level which are listed in table \ref{tab:tpi_M_K} for $K=1, ..., 10$. The interested reader is referred to \cite{Gear2003telescopic} for technical details on how to calculate these maximal values.

\begin{remark} \label{rem:rho_changing_M}
	Typically, we desire to fix the outermost time step $h_L=C\dx$ rather than choosing constant values $M$ on each level. This allows us to easily control the time instants at which the numerical solution is calculated. In that case, we first determine $L$ from equation \eqref{eq:tpi_zero_one_stable_L} using a fixed value of $K$ and selecting the corresponding maximal value of $M$ from table \ref{tab:tpi_M_K} as before. However, in general, for these values of $K$ and $M$ equation \eqref{eq:tpi_successive_timesteps} will not be equal to the chosen value $h_L$. To that end, we choose the value of $M_{L-1}$ on the outermost level as:
	\begin{equation} \label{eq:tpi_rho_M_outermost}
		M_{L-1} = \frac{h_L}{h_{L-1}} - K - 1, \qquad \mbox{if} \;\; \frac{h_L}{h_{L-1}} \ge K + 2,
	\end{equation}
	where the inequality on the right ensures that $M_{L-1} \ge 1$.
	However, when the chosen time step $h_L$ yields a value of $M_{L-1}$ in equation \eqref{eq:tpi_rho_M_outermost} less than 1 (or even negative), we decrease the value of $M$ on the lower levels, starting from the outermost to the innermost levels, until we find $M_{L-1} \ge 1$.
\end{remark}

\begin{table}
	\begin{center}
	    \begin{tabular}{ | c | c | c | c | c | c | c | c | c | c | c |}
	    \hline
	    $K$ & 1 & 2 & 3 & 4 & 5 & 6 & 7 & 8 & 9 & 10 \\ \hline
	    $M$ & 2 & 3 & 6.66 & 8.32 & 12.21 & 14.24 & 18.21 & 20.48 & 24.48 & 26.91 \\ \hline
	    \end{tabular}
	\end{center}
	\caption{\label{tab:tpi_M_K} Maximum value of $M$ for a given value of $K$ for a $[0,1]$-stable TPI method.}
\end{table}

When considering the time-varying spectrum case, we clearly find that both $M$ and $K$ are independent of $\epsi$. However, equation \eqref{eq:tpi_zero_one_stable_L} shows that the number of projective levels required for a $[0,1]$-stable TPI method increases as $O(\log(1/\epsi))$ given that $h_0=O(\epsi)$, see equation \eqref{eq:fe_h0_rho}. Therefore, the $[0,1]$-stable TPI method cost is not completely $\epsi$-independent but the dependence is rather modest.

In the numerical example, we choose a constant value $K=6$ on all projective levels. From table \ref{tab:tpi_M_K}, we deduce that the corresponding maximal value of $M$ on each projective level to obtain a $[0,1]$-stable TPFE method is 14.24. The required number of projective levels $L$ resulting from equation \eqref{eq:tpi_zero_one_stable_L} is 2. When choosing the outermost time step as $h_2 = 0.4\dx$, the adapted values of $M$ of the $[0,1]$-stable level-2 TPFE method are found as $M=\{14.24,11.79\}$, see remark \ref{rem:rho_changing_M}. The spectrum and stability region in each plane are shown in figure \ref{fig:spectrum_tpfe_rho}. Notice that, in every plane, the spectrum is not clustered anymore and is spread along the negative real axis. In addition, the stability region of the level-2 TPFE method does not split up into multiple disks such that it is indeed $[0,1]$-stable.

\begin{figure}[t]
	\begin{center}
		\figname{tprk4_spectrum_rho_4levels}
		\input{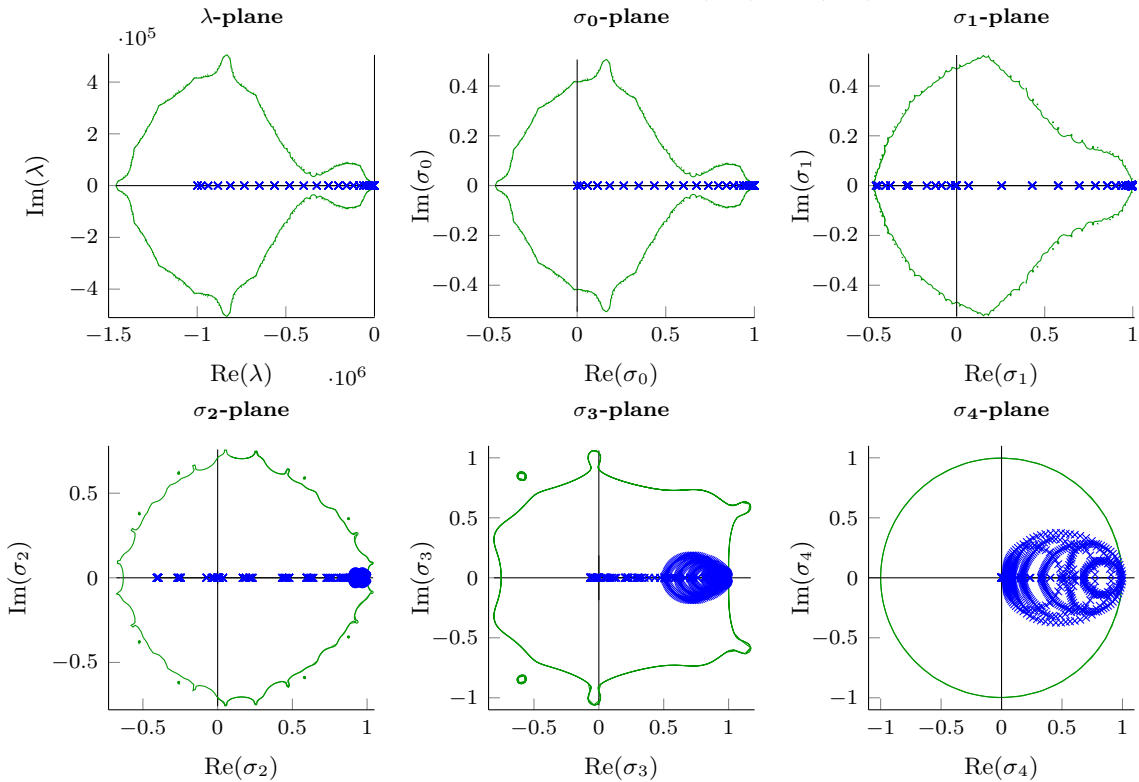}
	\end{center}
	\vspace{-0.4cm}\caption{\label{fig:spectrum_tprk4_rho} Stability analysis for a level-4 TPI method. Blue crosses in each plane correspond to eigenvalues in that plane. Green regions represent the stability regions of the full level-4 TPI method with respect to the plane under study. The red dashed regions are the stability regions of a classical projective integration method of the next level integrator in a certain plane. }
\end{figure}

Next, we repeat the same construction process for PRK4 as outermost integrator. As explained at the end of section \ref{subsec:tpi_params_omega}, the parameters of PRK4 are based on those of PFE as outermost integrator. We consider $\epsi=10^{-6}$ and fix $K=3$ on all projective levels. Then, the maximal value of $M$ is 6.66 and the required number of levels is $L=4$. When choosing the outermost time step as $h_4 = 0.4\dx$, the corrected values of $M$ of $[0,1]$-stable the level-4 TPRK4 method are $M=\{6.66, 6.26, 2.06, 2.03\}$. The results can be seen in figure \ref{fig:spectrum_tprk4_rho}. We conclude that the TPI construction procedure described in section \ref{subsec:tpi_params_rho} successfully results in a level-4 TPRK4 method for which the stability region does not split up.

\section{Numerical experiments} \label{sec:results}
We now examine the time stepping capabilities of the TPI method for equations of the form \eqref{eq:kin_eq_bgk} with $\collFreq(\x,t) = \rho(\x,t)$. We provide TPI construction and simulation test cases in 1D and 2D.

\paragraph{One-dimensional case ($\boldsymbol{D=1}$).}
First, we look at the one-dimensional kinetic equation \eqref{eq:kin_eq_bgk_1d} with linear Maxwellian given in equation \eqref{eq:maxwellian_artificial}. In that case, the limiting ${(\hydroLimit)}$ dynamics of equation \eqref{eq:kin_eq_bgk_1d} corresponds to the linear advection equation, given in equation \eqref{eq:lin_adv}. We compute the solution for $t \in [0,1]$ and $x \in [0,1]$. We impose periodic boundary conditions and choose a discontinuous initial density given by:
\begin{align} \label{eq:lin_adv_rho_init}
	\rho(x,0) = \left\{ \begin{array} {l@{\qquad \quad}l} 1 & 0.2 \le x < 0.4\\ 0.5 & 0.6 \le x < 0.8 \\ 0.1 & \text{otherwise}  \end{array} \right.&.
\end{align}
The initial distribution $f(x,v,0)$ is then chosen as the linearized Maxwellian given in equation \eqref{eq:maxwellian_artificial} corresponding to the initial density in \eqref{eq:lin_adv_rho_init}. We discretize velocity space using $J=10$ discrete velocity components obtained as the nodes of Gauss-Hermite quadrature with respect to the measure \eqref{eq:V_measure}. The innermost integrator is the forward Euler scheme with time step $h_0=\epsi$ and $\epsi=10^{-5}$. Since the initial density is sufficiently far from 0, we approximate the (linear) flux in equation \eqref{eq:kin_eq_bgk_1d} by the standard upwind differences of order 1, 2 and 3 with grid spacing $\dx=5 \cdot 10^{-3}$. We construct a $[0,1]$-stable TPRK4 method consisting of $L=2$ projective levels with constant $K=5$ and outermost time step $h_L = 0.5\dx$. The values of $M$ on each level are calculated as $M=\{12.21, 7.73\}$. 

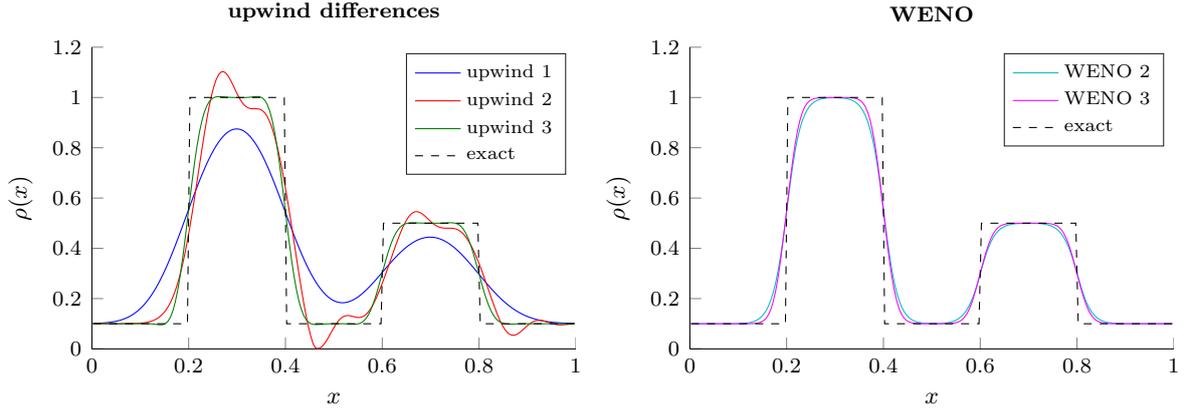
\begin{figure}[t]
	\begin{center}
		\figname{tpi_simulation_1d}
%
%
\definecolor{mycolor1}{rgb}{0.00000,0.49804,0.00000}%
\definecolor{mycolor2}{rgb}{0.00000,0.74902,0.74902}%
\definecolor{mycolor3}{rgb}{1.00000,0.00000,1.00000}%
\newcommand{\figWidth}{0.4\textwidth}
\newcommand{\figHeight}{4cm}
\newcommand{\figSpacingRight}{1.5cm} 
\begin{tikzpicture}

\begin{axis}[%
width=\figWidth,
height=\figHeight,
at={(0,0)},
scale only axis,
xmin=0,
xmax=1,
xlabel={$x$},
ymin=0,
ymax=1.2,
ytick={0, 0.2, ..., 1.2},
ylabel={$\rho(x)$},
axis background/.style={fill=white},
title style={font=\bfseries\scriptsize},
title={upwind differences},
axis x line*=bottom,
axis y line*=left,
legend style={legend cell align=left,align=left,draw=white!15!black}
]
\addplot [color=blue,solid]
  table[]{tikz/data/simulation_1d-1.tsv};
\addlegendentry{upwind 1};

\addplot [color=red,solid]
  table[]{tikz/data/simulation_1d-2.tsv};
\addlegendentry{upwind 2};

\addplot [color=mycolor1,solid]
  table[]{tikz/data/simulation_1d-3.tsv};
\addlegendentry{upwind 3};

\addplot [color=black,dashed]
  table[]{tikz/data/simulation_1d-4.tsv};
\addlegendentry{exact};

\end{axis}

\begin{axis}[%
width=\figWidth,
height=\figHeight,
at={(\figWidth+\figSpacingRight,0)},
scale only axis,
xmin=0,
xmax=1,
xlabel={$x$},
ymin=0,
ymax=1.2,
ytick={0, 0.2, ..., 1.2},
ylabel={$\rho(x)$},
axis background/.style={fill=white},
title style={font=\bfseries\scriptsize},
title={WENO},
axis x line*=bottom,
axis y line*=left,
legend style={legend cell align=left,align=left,draw=white!15!black}
]
\addplot [color=mycolor2,solid]
  table[]{tikz/data/simulation_1d-5.tsv};
\addlegendentry{WENO 2};

\addplot [color=mycolor3,solid]
  table[]{tikz/data/simulation_1d-6.tsv};
\addlegendentry{WENO 3};

\addplot [color=black,dashed]
  table[]{tikz/data/simulation_1d-7.tsv};
\addlegendentry{exact};

\end{axis}
\end{tikzpicture}%
	\end{center}
	\vspace{-0.4cm}\caption{\label{fig:tpi_sim_1d} Simulation results (density) at $t=1$ of equation \eqref{eq:kin_eq_bgk_1d} with $\collFreq(x,t)=\rho(x,t)$ using a level-2 TPRK4 method with time step $h_2 = 0.5\dx$ and $\dx = 5 \cdot 10^{-3}$. We compare two different spatial discretization techniques. Left: upwind differences of order 1 (blue), 2 (red) and 3 (green). Right: WENO2 (cyan) and WENO3 (purple). The black dashed curve in both plots corresponds to the exact solution of the linear advection equation \eqref{eq:lin_adv}. }
\end{figure}

The numerical solution at $t=1$ for different orders of the upwind scheme is shown in the left plot of figure \ref{fig:tpi_sim_1d}. We also plotted the exact solution of the limiting linear advection equation in black. As can be seen, the first order upwind method is too diffusive and its higher-order versions produce spurious oscillations around discontinuities. To counter this undesired result, we implemented a Weighted Essentially Non-Oscillatory (WENO) scheme \cite{eno_weno} which uses a weighted linear combination of all possible stencils for each grid point for a given spatial order of accuracy giving more weight to smooth stencils. The results for WENO2 and WENO3 are depicted by cyan and purple lines, respectively, on the right plot of figure \ref{fig:tpi_sim_1d}. In this case, we obtain a high order approximation without oscillations.

\begin{figure}[t]
	\begin{center}
		\figname{tpi_simulation_2d}
%
%
\newcommand{\figWidth}{3.7cm} 
\newcommand{\figHeight}{\figWidth} 
\newcommand{\figSpacingRight}{1.5cm} 
\begin{tikzpicture}

\begin{axis}[%
width=\figWidth,
height=\figHeight,
at={(0,0.0)},
scale only axis,
plot box ratio=1 1 1,
xmin=0,
xmax=1.01,
tick align=outside,
xlabel={$x$},
xmajorgrids,
ymin=0,
ymax=1,
ylabel={$y$},
ymajorgrids,
zmin=-0.01,
zmax=1,
zlabel={$\rho(\x)$},
zmajorgrids,
view={-45}{0},
axis background/.style={fill=white},
title style={font=\scriptsize\bfseries},
title={upwind differences order 1},
axis x line*=bottom,
axis y line*=left,
axis z line*=left
]

\addplot3[%
surf,
shader=flat corner,draw=black,opacity=0.3,fill opacity=1,z buffer=sort,colormap={mymap}{[1pt] rgb(0pt)=(0.2081,0.1663,0.5292); rgb(1pt)=(0.211624,0.189781,0.577676); rgb(2pt)=(0.212252,0.213771,0.626971); rgb(3pt)=(0.2081,0.2386,0.677086); rgb(4pt)=(0.195905,0.264457,0.7279); rgb(5pt)=(0.170729,0.291938,0.779248); rgb(6pt)=(0.125271,0.324243,0.830271); rgb(7pt)=(0.0591333,0.359833,0.868333); rgb(8pt)=(0.0116952,0.38751,0.881957); rgb(9pt)=(0.00595714,0.408614,0.882843); rgb(10pt)=(0.0165143,0.4266,0.878633); rgb(11pt)=(0.0328524,0.443043,0.871957); rgb(12pt)=(0.0498143,0.458571,0.864057); rgb(13pt)=(0.0629333,0.47369,0.855438); rgb(14pt)=(0.0722667,0.488667,0.8467); rgb(15pt)=(0.0779429,0.503986,0.838371); rgb(16pt)=(0.0793476,0.520024,0.831181); rgb(17pt)=(0.0749429,0.537543,0.826271); rgb(18pt)=(0.0640571,0.556986,0.823957); rgb(19pt)=(0.0487714,0.577224,0.822829); rgb(20pt)=(0.0343429,0.596581,0.819852); rgb(21pt)=(0.0265,0.6137,0.8135); rgb(22pt)=(0.0238905,0.628662,0.803762); rgb(23pt)=(0.0230905,0.641786,0.791267); rgb(24pt)=(0.0227714,0.653486,0.776757); rgb(25pt)=(0.0266619,0.664195,0.760719); rgb(26pt)=(0.0383714,0.674271,0.743552); rgb(27pt)=(0.0589714,0.683757,0.725386); rgb(28pt)=(0.0843,0.692833,0.706167); rgb(29pt)=(0.113295,0.7015,0.685857); rgb(30pt)=(0.145271,0.709757,0.664629); rgb(31pt)=(0.180133,0.717657,0.642433); rgb(32pt)=(0.217829,0.725043,0.619262); rgb(33pt)=(0.258643,0.731714,0.595429); rgb(34pt)=(0.302171,0.737605,0.571186); rgb(35pt)=(0.348167,0.742433,0.547267); rgb(36pt)=(0.395257,0.7459,0.524443); rgb(37pt)=(0.44201,0.748081,0.503314); rgb(38pt)=(0.487124,0.749062,0.483976); rgb(39pt)=(0.530029,0.749114,0.466114); rgb(40pt)=(0.570857,0.748519,0.44939); rgb(41pt)=(0.609852,0.747314,0.433686); rgb(42pt)=(0.6473,0.7456,0.4188); rgb(43pt)=(0.683419,0.743476,0.404433); rgb(44pt)=(0.71841,0.741133,0.390476); rgb(45pt)=(0.752486,0.7384,0.376814); rgb(46pt)=(0.785843,0.735567,0.363271); rgb(47pt)=(0.818505,0.732733,0.34979); rgb(48pt)=(0.850657,0.7299,0.336029); rgb(49pt)=(0.882433,0.727433,0.3217); rgb(50pt)=(0.913933,0.725786,0.306276); rgb(51pt)=(0.944957,0.726114,0.288643); rgb(52pt)=(0.973895,0.731395,0.266648); rgb(53pt)=(0.993771,0.745457,0.240348); rgb(54pt)=(0.999043,0.765314,0.216414); rgb(55pt)=(0.995533,0.786057,0.196652); rgb(56pt)=(0.988,0.8066,0.179367); rgb(57pt)=(0.978857,0.827143,0.163314); rgb(58pt)=(0.9697,0.848138,0.147452); rgb(59pt)=(0.962586,0.870514,0.1309); rgb(60pt)=(0.958871,0.8949,0.113243); rgb(61pt)=(0.959824,0.921833,0.0948381); rgb(62pt)=(0.9661,0.951443,0.0755333); rgb(63pt)=(0.9763,0.9831,0.0538)},mesh/rows=25]
table[point meta=\thisrow{c}] {%
tikz/data/simulation_2d-1.tsv};
\end{axis}

\begin{axis}[%
width=\figWidth,
height=\figHeight,
at={(\figWidth+\figSpacingRight,0)},
scale only axis,
plot box ratio=1 1 1,
xmin=0,
xmax=1.01,
tick align=outside,
xlabel={$x$},
xmajorgrids,
ymin=0,
ymax=1,
ylabel={$y$},
ymajorgrids,
zmin=-0.01,
zmax=1,
zlabel={$\rho(\x)$},
zmajorgrids,
view={-45}{0},
axis background/.style={fill=white},
title style={font=\scriptsize\bfseries},
title={WENO2},
axis x line*=bottom,
axis y line*=left,
axis z line*=left
]

\addplot3[%
surf,
shader=flat corner,draw=black,opacity=0.3,fill opacity=1,z buffer=sort,colormap={mymap}{[1pt] rgb(0pt)=(0.2081,0.1663,0.5292); rgb(1pt)=(0.211624,0.189781,0.577676); rgb(2pt)=(0.212252,0.213771,0.626971); rgb(3pt)=(0.2081,0.2386,0.677086); rgb(4pt)=(0.195905,0.264457,0.7279); rgb(5pt)=(0.170729,0.291938,0.779248); rgb(6pt)=(0.125271,0.324243,0.830271); rgb(7pt)=(0.0591333,0.359833,0.868333); rgb(8pt)=(0.0116952,0.38751,0.881957); rgb(9pt)=(0.00595714,0.408614,0.882843); rgb(10pt)=(0.0165143,0.4266,0.878633); rgb(11pt)=(0.0328524,0.443043,0.871957); rgb(12pt)=(0.0498143,0.458571,0.864057); rgb(13pt)=(0.0629333,0.47369,0.855438); rgb(14pt)=(0.0722667,0.488667,0.8467); rgb(15pt)=(0.0779429,0.503986,0.838371); rgb(16pt)=(0.0793476,0.520024,0.831181); rgb(17pt)=(0.0749429,0.537543,0.826271); rgb(18pt)=(0.0640571,0.556986,0.823957); rgb(19pt)=(0.0487714,0.577224,0.822829); rgb(20pt)=(0.0343429,0.596581,0.819852); rgb(21pt)=(0.0265,0.6137,0.8135); rgb(22pt)=(0.0238905,0.628662,0.803762); rgb(23pt)=(0.0230905,0.641786,0.791267); rgb(24pt)=(0.0227714,0.653486,0.776757); rgb(25pt)=(0.0266619,0.664195,0.760719); rgb(26pt)=(0.0383714,0.674271,0.743552); rgb(27pt)=(0.0589714,0.683757,0.725386); rgb(28pt)=(0.0843,0.692833,0.706167); rgb(29pt)=(0.113295,0.7015,0.685857); rgb(30pt)=(0.145271,0.709757,0.664629); rgb(31pt)=(0.180133,0.717657,0.642433); rgb(32pt)=(0.217829,0.725043,0.619262); rgb(33pt)=(0.258643,0.731714,0.595429); rgb(34pt)=(0.302171,0.737605,0.571186); rgb(35pt)=(0.348167,0.742433,0.547267); rgb(36pt)=(0.395257,0.7459,0.524443); rgb(37pt)=(0.44201,0.748081,0.503314); rgb(38pt)=(0.487124,0.749062,0.483976); rgb(39pt)=(0.530029,0.749114,0.466114); rgb(40pt)=(0.570857,0.748519,0.44939); rgb(41pt)=(0.609852,0.747314,0.433686); rgb(42pt)=(0.6473,0.7456,0.4188); rgb(43pt)=(0.683419,0.743476,0.404433); rgb(44pt)=(0.71841,0.741133,0.390476); rgb(45pt)=(0.752486,0.7384,0.376814); rgb(46pt)=(0.785843,0.735567,0.363271); rgb(47pt)=(0.818505,0.732733,0.34979); rgb(48pt)=(0.850657,0.7299,0.336029); rgb(49pt)=(0.882433,0.727433,0.3217); rgb(50pt)=(0.913933,0.725786,0.306276); rgb(51pt)=(0.944957,0.726114,0.288643); rgb(52pt)=(0.973895,0.731395,0.266648); rgb(53pt)=(0.993771,0.745457,0.240348); rgb(54pt)=(0.999043,0.765314,0.216414); rgb(55pt)=(0.995533,0.786057,0.196652); rgb(56pt)=(0.988,0.8066,0.179367); rgb(57pt)=(0.978857,0.827143,0.163314); rgb(58pt)=(0.9697,0.848138,0.147452); rgb(59pt)=(0.962586,0.870514,0.1309); rgb(60pt)=(0.958871,0.8949,0.113243); rgb(61pt)=(0.959824,0.921833,0.0948381); rgb(62pt)=(0.9661,0.951443,0.0755333); rgb(63pt)=(0.9763,0.9831,0.0538)},mesh/rows=25]
table[point meta=\thisrow{c}] {%
tikz/data/simulation_2d-2.tsv};
\end{axis}

\begin{axis}[%
width=\figWidth,
height=\figHeight,
at={(2*\figWidth+2*\figSpacingRight,0)},
scale only axis,
plot box ratio=1 1 1,
xmin=0,
xmax=1.01,
tick align=outside,
xlabel={$x$},
xmajorgrids,
ymin=0,
ymax=1,
ylabel={$y$},
ymajorgrids,
zmin=-0.01,
zmax=1,
zlabel={$\rho(\x)$},
zmajorgrids,
view={-45}{0},
axis background/.style={fill=white},
title style={font=\scriptsize\bfseries},
title={WENO3},
axis x line*=bottom,
axis y line*=left,
axis z line*=left
]

\addplot3[%
surf,
shader=flat corner,draw=black,opacity=0.3,fill opacity=1,z buffer=sort,colormap={mymap}{[1pt] rgb(0pt)=(0.2081,0.1663,0.5292); rgb(1pt)=(0.211624,0.189781,0.577676); rgb(2pt)=(0.212252,0.213771,0.626971); rgb(3pt)=(0.2081,0.2386,0.677086); rgb(4pt)=(0.195905,0.264457,0.7279); rgb(5pt)=(0.170729,0.291938,0.779248); rgb(6pt)=(0.125271,0.324243,0.830271); rgb(7pt)=(0.0591333,0.359833,0.868333); rgb(8pt)=(0.0116952,0.38751,0.881957); rgb(9pt)=(0.00595714,0.408614,0.882843); rgb(10pt)=(0.0165143,0.4266,0.878633); rgb(11pt)=(0.0328524,0.443043,0.871957); rgb(12pt)=(0.0498143,0.458571,0.864057); rgb(13pt)=(0.0629333,0.47369,0.855438); rgb(14pt)=(0.0722667,0.488667,0.8467); rgb(15pt)=(0.0779429,0.503986,0.838371); rgb(16pt)=(0.0793476,0.520024,0.831181); rgb(17pt)=(0.0749429,0.537543,0.826271); rgb(18pt)=(0.0640571,0.556986,0.823957); rgb(19pt)=(0.0487714,0.577224,0.822829); rgb(20pt)=(0.0343429,0.596581,0.819852); rgb(21pt)=(0.0265,0.6137,0.8135); rgb(22pt)=(0.0238905,0.628662,0.803762); rgb(23pt)=(0.0230905,0.641786,0.791267); rgb(24pt)=(0.0227714,0.653486,0.776757); rgb(25pt)=(0.0266619,0.664195,0.760719); rgb(26pt)=(0.0383714,0.674271,0.743552); rgb(27pt)=(0.0589714,0.683757,0.725386); rgb(28pt)=(0.0843,0.692833,0.706167); rgb(29pt)=(0.113295,0.7015,0.685857); rgb(30pt)=(0.145271,0.709757,0.664629); rgb(31pt)=(0.180133,0.717657,0.642433); rgb(32pt)=(0.217829,0.725043,0.619262); rgb(33pt)=(0.258643,0.731714,0.595429); rgb(34pt)=(0.302171,0.737605,0.571186); rgb(35pt)=(0.348167,0.742433,0.547267); rgb(36pt)=(0.395257,0.7459,0.524443); rgb(37pt)=(0.44201,0.748081,0.503314); rgb(38pt)=(0.487124,0.749062,0.483976); rgb(39pt)=(0.530029,0.749114,0.466114); rgb(40pt)=(0.570857,0.748519,0.44939); rgb(41pt)=(0.609852,0.747314,0.433686); rgb(42pt)=(0.6473,0.7456,0.4188); rgb(43pt)=(0.683419,0.743476,0.404433); rgb(44pt)=(0.71841,0.741133,0.390476); rgb(45pt)=(0.752486,0.7384,0.376814); rgb(46pt)=(0.785843,0.735567,0.363271); rgb(47pt)=(0.818505,0.732733,0.34979); rgb(48pt)=(0.850657,0.7299,0.336029); rgb(49pt)=(0.882433,0.727433,0.3217); rgb(50pt)=(0.913933,0.725786,0.306276); rgb(51pt)=(0.944957,0.726114,0.288643); rgb(52pt)=(0.973895,0.731395,0.266648); rgb(53pt)=(0.993771,0.745457,0.240348); rgb(54pt)=(0.999043,0.765314,0.216414); rgb(55pt)=(0.995533,0.786057,0.196652); rgb(56pt)=(0.988,0.8066,0.179367); rgb(57pt)=(0.978857,0.827143,0.163314); rgb(58pt)=(0.9697,0.848138,0.147452); rgb(59pt)=(0.962586,0.870514,0.1309); rgb(60pt)=(0.958871,0.8949,0.113243); rgb(61pt)=(0.959824,0.921833,0.0948381); rgb(62pt)=(0.9661,0.951443,0.0755333); rgb(63pt)=(0.9763,0.9831,0.0538)},mesh/rows=25]
table[point meta=\thisrow{c}] {%
tikz/data/simulation_2d-3.tsv};
\end{axis}
\end{tikzpicture}%
	\end{center}
	\vspace{-0.4cm}\caption{\label{fig:tpi_sim_2d} Simulation results (density) at $t=1$ of equation \eqref{eq:kin_eq_bgk_2d} with $\collFreq(\x,t)=\rho(\x,t)$ using a level-3 TPRK4 method with time step $h_2 = 0.5\dx$ and $\dx = \dy = 0.02$. We compare three spatial discretization techniques: first-order upwind (left), WENO2 (middle) and WENO3 (right). }
\end{figure}
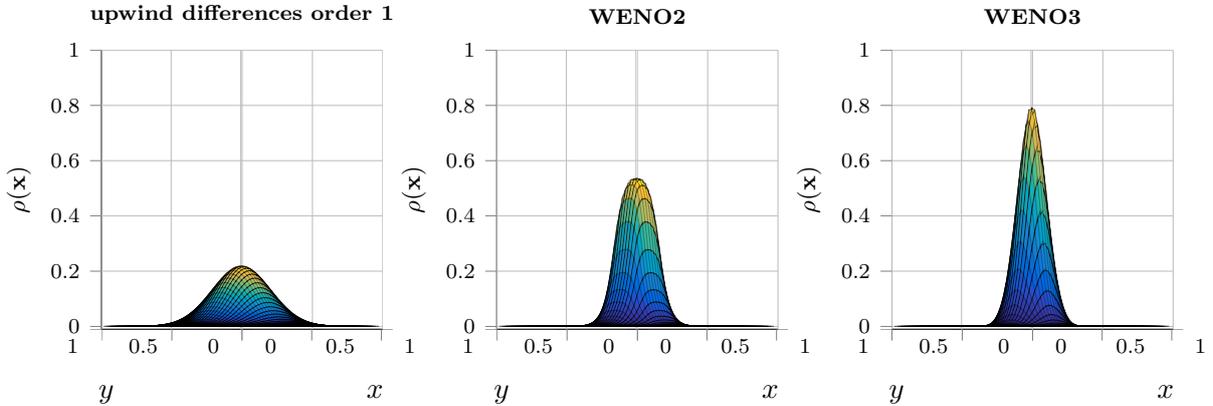

\paragraph{Two-dimensional case ($\boldsymbol{D=2}$).}
Next, we examine the two-dimensional equation \eqref{eq:kin_eq_bgk_2d} with linear Maxwellian given in \eqref{eq:maxwellian_artificial_2d}. Now, the limiting $(\hydroLimit)$ dynamics corresponds to the two-dimensional linear advection equation as shown in \eqref{eq:lin_adv_2d}. We again compute the solution for $t \in [0,1]$ and consider a rectangular spatial domain $(x,y) \in [0,1] \times [0,1]$. We impose periodic boundary conditions and choose a smooth initial density given by a domain-centered Gaussian function:
\begin{equation} \label{eq:lin_adv_rho_init_2d}
	\rho(\x,0) = \exp\left( -100\abs{\x - 0.5}^2 \right).
\end{equation}
We discretize velocity space $(v^x,v^y) \subset \R^2$ using $10 \times 10$ discrete velocity components obtained as the nodes of Gauss-Hermite quadrature with respect to the two-dimensional measure \eqref{eq:V_measure_2d}. In the TPI framework, we choose the forward Euler scheme with time step $h_0=\epsi$ and $\epsi=10^{-5}$ as innermost integrator. The spatial domain is discretized by choosing a rectangular grid with grid spacing $\dx = \dy = 0.02$. Then, the (linear) fluxes in equation \eqref{eq:kin_eq_bgk_2d} are approximated by first order upwind differences. We construct a $[0,1]$-stable TPRK4 method consisting of $L=3$ projective levels with constant $K=3$ and outermost time step $h_L = 0.5\dx$. The values of $M$ on each level are calculated as $M=\{6.66, 6.66, 4.81\}$. The result is shown in the left plot of figure \ref{fig:tpi_sim_2d}.

Clearly, the first-order upwind method again introduces very strong numerical diffusion. Therefore, we turn towards higher-order spatial discretization techniques. The obvious choice of higher-order upwind methods is inappropriate for the considered initial solution, since these inevitably generate under- and overshoots in the numerical solution causing a potential loss of stability in finite time, see equation \eqref{eq:maximum_principle}. However, we can obtain higher-order solutions by applying the WENO scheme. The result for WENO2 and WENO3 is visualized in the middle and right plots, respectively, of figure \ref{fig:tpi_sim_2d}.

\section{Conclusions} \label{sec:conclusions}
We presented a general, higher-order, fully explicit integration method for kinetic equations with BGK-like source term containing a collision frequency leading to multiple relaxation times. The method uses a hierarchy of projective integrators ,leading to telescopic projective integration methods.  The number of levels, as well as the size and number of the time-steps at each level, can be derived based on the problem's spectrum. Unlike other methods based on relaxation \cite{Jin1995,Aregba-Driollet2000}, the telescopic projective integration method does not rely on a splitting technique, but only on an appropriate selection of time steps using a naive explicit discretization method at its core. Its main advantage is its generality and ease of use.

We showed that, with an appropriate choice of the inner time step, the time step restriction on the outer time step is independent of the small-scale parameter. Moreover, the number of inner integrator iterations and the projective step size are independent of the scaling parameter. By contrast, the required number of projective levels depends on this parameter, be it only logarithmically. We analyzed stability and provided numerical results on the method construction procedure. We applied the method both to one- and two-dimensional kinetic equations.

In future work, we foresee to construct stable telescopic projective integration methods for the nonlinear BGK kinetic equation without linearizing the Maxwellian distribution. 
Following that, an ambitious goal is to apply these methods to the full Boltzmann equation, for which the spectrum of the corresponding linearized collision operator is known to consist of a number of separated eigenvalue clusters, see \cite{Ellis1975,Rey2013}. However, a more precise characterization of the location and size of these clusters is required to determine suitable method parameters.


\section*{Acknowledgement}
We would like to thank Thomas Rey from Laboratoire Paul Painlev\'e of Universit\'e de Lille for providing us with clear background material on Boltzmann and BGK kinetic equations and his assistance with required derivations and calculations which supplemented the motivation of this work.

\bibliographystyle{plain}

\bibliography{refs_mendeley}

\end{document}